\documentclass[11pt]{article}
\usepackage{graphicx} 
\usepackage{geometry}
\usepackage{amsmath}
\usepackage{amssymb}
\usepackage{amsfonts}
\usepackage{multirow}
\usepackage{amsthm}
\usepackage{xcolor}
\usepackage[numbers]{natbib}
\usepackage{comment}
\usepackage{enumitem}
\usepackage{hyperref}
\usepackage{algorithm}
\usepackage{algpseudocode}
\usepackage{rotating}
\usepackage{centernot}
\newtheorem{theorem}{Theorem}
\newtheorem{lemma}{Lemma}
\newtheorem{corollary}{Corollary}
\newtheorem{proposition}{Proposition}
\theoremstyle{definition}
\newtheorem{definition}{Definition}
\newtheorem{assumption}{A}
\newtheorem{example}{Example}

 {\begin{list}{}%
         {\setlength{\leftmargin}{#1}}%
         \item[]%
 }
 {\end{list}}

\newcommand{\E}{\text{E}} 

\newcommand{\di}{{\rm d}}

\newcommand{\N}{\mbox{N}}

\newcommand{\cF}{{\cal F}}

\newcommand{\seq}[1]{({#1})_{n \geq 0}}
\graphicspath{{figures/}}

\newcommand{\cG}{{\cal G}}
\newcommand{\nll}{\centernot{\ll}}

\newif\ifrevision
\revisiontrue

\ifrevision
    
    \newcommand{\removed}[1]{\sout{#1}}
    \newcommand{\commentforreviewer}[1]{\textcolor{blue}{#1}}
\else
    
    \newcommand{\removed}[1]{}
    \newcommand{\commentforreviewer}[1]{}
\fi

\ExplSyntaxOn
\NewDocumentEnvironment{multiequation}{b}
 {
  \vantiempham:n { #1 }
 }
 {}
\seq_new:N \l__vantiempham_md_seq
\cs_new_protected:Nn \vantiempham:n
 {
  \seq_set_split:Nnn \l__vantiempham_md_seq { \\ } { #1 }
  $$
  \seq_map_function:NN \l__vantiempham_md_seq \__vantiempham_md_item:n
  $$
 }
\cs_new_protected:Nn \__vantiempham_md_item:n
 {
  \begin{minipage}{\dim_eval:n {\displaywidth/(\seq_count:N \l__vantiempham_md_seq)} }
  \begin{equation}
  \cs_set_eq:Nc \label { ltx@label }
  #1
  \vphantom{\seq_use:Nn \l__vantiempham_md_seq {}}
  \end{equation}
  \end{minipage}
 }
\ExplSyntaxOff

\title{Limit theorems for almost conditionally identically distributed random variables} %
 \author{Marco Battiston$^{1}$ and Lorenzo Cappello$^{2,3}$\\
 		$^1$School of Mathematical Sciences, Lancaster University\\
	$^2$ Department of Economics and Business, Universitat Pompeu Fabra\\
$^3$ Data Science Center, Barcelona School of Economics}

\begin{document}

\maketitle.                      
\begin{abstract}
$\xi$-Almost Conditional Identically Distributed ($\xi$-a.c.i.d.) random variables are generalizations of conditional identically distributed (c.i.d.) and exchangeable ones. $\xi$-a.c.i.d. random variables arise naturally in applications in statistics, such as in recursive algorithms 
 and heteroskedastic observations, and are parametrized by a sequence $\xi=(\xi_{n})_{n\geq 0}$ that quantifies departure 
from the c.i.d. property. In this work, we extend limit theorems for exchangeable and c.i.d. random variables to the more general setting of $\xi$-a.c.i.d. random variables. After proving a Strong Law of Large Numbers, we provide necessary and sufficient conditions for the asymptotic directing measure of the sequence to be absolutely continuous with respect to a given $\sigma$-finite measure. Then, under regularity assumptions on the parameter $\xi$, three different Central Limit Theorems for $\xi$-a.c.i.d. sequences are derived, involving the predictive, empirical and asymptotic directing distributions. These theorems have several statistical applications, including the construction of asymptotic credible intervals in Bayesian predictive inference. 
\end{abstract}



\section{Introduction}

Let $X=(X_{n})_{n\geq 1}$ be a sequence of random variables on a Polish space $(E,\mathcal{E})$, $\mathcal{G}=(\mathcal{G}_{n})_{n\geq 0}$ its natural filtration, with $\mathcal{G}_{0}=\{\emptyset, E\}$, and $\xi=(\xi_{n})_{n\geq 0}$ a sequence of random variables  adapted to $\mathcal{G}$. We define $X$ to be $\xi$-a.c.i.d. if the sequence of predictive distributions $(\mathbb{P}(X_{n+1}\in \cdot |\mathcal{G}_{n}))_{n\geq 0}$ forms a measure-valued almost supermartingale with parameters $\xi=(\xi_{n})_{n\geq 0}$, see Definition~\ref{def:acid.pred}.  As an important particular case, when $\xi_{n}=0$ a.s. for all $n\geq 0$, the sequence of predictive distribution is a martingale, and the sequence $X$ is conditionally identically distributed (c.i.d.) \citep{kal88,berti2004limit}. If, in addition, the sequence is also stationary, then it is exchangeable \citep{kal88}. 

C.i.d. random variables have been extensively studied in the probability and statistics literature, see e.g. \cite{berti2004limit,ber12,bcl10,bw25,fps18,fp25,bdpr21,bdlpr25}, and  have gained attention very recently in the context of Bayesian predictive inference \citep{fong2021martingale,fortini2020quasi,bdpr21}; see \citep{bdlpr25,fp25} for a review. Within this approach, a statistician bypasses both the need to specify a a sampling model and a prior for the unknown parameter, which might be complicated in complex data, and to use Markov Chain Monte Carlo for posterior inference, which is typically computationally intensive. Instead, the practitioner needs to specify only a sequence of predictive distributions $\left(\mathbb{P}(X_{n+1}\in \cdot | \mathcal{G}_{n})\right)_{n\geq 0}$. Inference can be then performed via a simulation based scheme called \emph{Predictive Resampling}, a Bootstrap-like algorithm; see \cite{fong2021martingale,fortini2020quasi}. Alternatively, when applicable, limit theorems can be 
used to derive asymptotic approximations to credible intervals; see \cite{fortini2020quasi, fortini2023prediction}. A key question is how to specify $\left(\mathbb{P}(X_{n+1}\in \cdot | \mathcal{G}_{n})\right)_{n\geq 0}$ without resorting to the usual prior/posterior scheme. Much work has focused on fast recursive rules, and the literature is rapidly growing. To mention a few examples, some works have studied conditions under which the corresponding sequences of observables remain exchangeable \citep{berti2023kernel,sariev2024characterization,sariev2025sufficientness,chorbadzhiyska2025some}, while others have considered the broader c.i.d. class \cite{bdpr21,bw25,fong2021martingale,fortini2020quasi}. A common thread is that these works  largely rely on limit theorems for exchangeable random variables (e.g., \citep{ald85,kal88,for00}) and c.i.d. ones \citep{berti2004limit,bpr13}.

Recent contributions have considered predictive schemes that induce laws for the observables that are not necessarily c.i.d.; see, for example, \cite{holmes2023statistical,fong2026asymptotics,bc25,lpr26,garelli2024asymptotics,glpr26b}. This calls for tools that go beyond the classical c.i.d. framework, and \cite{bc25,lpr26} have taken steps in this direction. In particular, \cite{bc25} show that several recently proposed predictive schemes generate $\xi$-a.c.i.d. sequences. This highlights the importance of extending limit theorems to this more flexible class of random variables. Examples of $\xi$-a.c.i.d. sequences include predictive distributions based kernel density estimators \citep{bc25} and the parametric Bayesian bootstrap \citep{holmes2023statistical,fong2026asymptotics}. We next provide a further example of a $\xi$-a.c.i.d. sequence, whose $\xi$-a.c.i.d. property has not previously been established in the literature:

\begin{example}\label{ex:garelli}
\cite{garelli2024asymptotics} proposes a general sequence of predictive distributions driven by sample mean and variance estimators. Here, we consider the specific instance in which the predictive distributions are Gaussian. Specifically, the marginal of the first observation is $X_1 \sim N(M_0, Q_0)$ on $\mathbb{R}^p$, where $M_{0}$ and $Q_{0}$ are some initial fixed mean vector and covariance matrix. Then, \cite{garelli2024asymptotics} defines recursively the predictive distributions, for $n \ge 1$, as follows
\[
X_{n+1} \mid X_{1:n} \sim N(M_n, Q_n), \qquad M_n = \frac{1}{n}\sum_{i=1}^n X_i, \qquad
Q_n = \frac{1}{n}\sum_{i=1}^n (X_i - M_n)(X_i - M_n)^T.
\]    
\end{example}
This sequence of predictive distributions can be shown to be $\xi$-a.c.i.d., with $\xi_{n}:=\frac{C_{p}}{(n+1)^{2}}$, where $C_{p}$ is a known constant depending only on $p$. The proof is in the Appendix~\ref{sec:proof.garelli}.

 A sufficient condition for a sequence of predictive distributions $\left(\mathbb{P}(X_{n+1}\in \cdot | \mathcal{G}_{n})\right)_{n\geq 0}$ to be suitable for Bayesian predictive inference is that it corresponds to an \emph{asymptotically exchangeable} sequence. Heuristically, this means that after discarding some initial observations, i.e. repeatedly applying the right shift transformation to the sequence, the process converges in distribution to an exchangeable sequence, see \cite{ald85} 
 for a formal definition. Informally, with some abuse of terminology, we also say that $(X_n)_{n\geq1}$ has the same limit random directing measure, denoted  $\tilde{F}$, as the limiting exchangeable sequence. We refer to \citep{fp25} for further details on why this property allows one to attach a Bayesian interpretation to $\tilde{F}$ and its distribution. In Theorem~\ref{thm:asym.exch}, we show that any $\xi$-a.c.i.d. sequence with a.s. summable $(\xi_{n})_{n\geq 0}$ is asymptotically exchangeable, hence amenable to Bayesian predictive inference. 

The goal of this paper is to provide limit theorems 
for generic $\xi$-a.c.i.d. sequences, extending corresponding results available for exchangeable \citep{ald85}, and c.i.d. random variables \citep{berti2004limit,bpr13}, respectively. Specifically, we firstly derive a Strong Law of Large Numbers (SLLN) and a new proof of asymptotic exchangeability for $\xi$-a.c.i.d. sequences. Secondly, Theorem~\ref{thm:abs.cont} derives necessary and sufficient conditions for the random measure $\tilde{F}$ to be a.s. absolutely continuous with respect to a given $\sigma$-finite measure $\lambda$. Some interesting connections to Theorem 1 of \cite{bpr13} are highlighted. 
Thirdly, three distinct Central Limit Theorems (CLTs) based on the predictive and empirical distributions of a generic $\xi$-a.c.i.d. sequence are derived. The proofs of these theorems rely on generalizations of classic martingale CLTs, which are covered in the Appendix. The relevance of these CLTs to applications, particularly Bayesian predictive inference and the construction of random intervals that can be interpreted as Bayesian credible intervals, is also discussed.

After a summary of all the notation in the next paragraph, the rest of the paper is organized as follows.  In Section~\ref{sec:acid}, we review the definitions of $\xi$-a.c.i.d. random variables and present some examples.
In Section~\ref{sec:slln.ae}, we prove a SLLN and asymptotic exchangeability. Section~\ref{sec:abs.cont} deals with the absolute continuity of the asymptotic directing measure $\tilde{F}$. Section~\ref{sec:CLTs} discusses three CLTs for $\xi$-a.c.i.d. sequences. Finally, some background and auxiliary results on almost supermartingales and CLTs are postponed to the Appendix. \vspace{0.3cm}

\noindent \textbf{Notation}: All random variables (r.v.s) are assumed to be defined on a common probability space $(\Omega,\mathcal{A},\mathbb{P})$.
$(X_{n})_{n\geq 1}$ denotes a sequence of r.v.s, and $\mathcal{G}=(\mathcal{G}_{n})_{n\geq 0}$, with $\mathcal{G}_{0}=\{\emptyset,\Omega\}$, its natural filtration. 
$\{\mathcal{F}_{n,i};n\geq 0, 0\leq i \leq k_{n} \}$ denotes a (possibly triangular) array of sub-$\sigma$-fields of $\mathcal{A}$, and $(Y_{n,i})_{n\geq 0,0\leq i \leq k_{n}}$  a (possibly triangular) array of r.v.s, with $Y_{n,i}$ assumed $\mathcal{F}_{n,i}$ measurable.
All $X_{n}$'s and $Y_{n,i}$'s take values on a Polish space $(E,\mathcal{E})$. $X_{1:n}$ is used as shorthand notation for the random vector $(X_{1},\ldots,X_{n})$. For a random variable or vector $X$, $\mathcal{L}(X)$ denotes the probability law of $X$. 
$\stackrel{\text{w}}{\to}$, $\stackrel{\text{d}}{\to}$ and $\stackrel{\text{a.s.}}{\to}$ 
stand for weak convergence of probability measures, convergence in distribution and almost sure convergence of r.v.s, respectively. $\mathbb{I}(\cdot)$ denotes the indicator function of the event within brackets, while $\delta_{X}(\cdot)$ denotes the Delta Dirac measure at $X$. For a given measure $\lambda$, $\lambda^{n}$ denotes its $n$-fold product measure.
$P_{n}(\omega,A) := \mathbb{P}\left( X_{n+1}\in A | \mathcal{G}_{n} \right)(\omega)$ denotes a regular version of  the conditional distribution of $X_{n+1}$ given $\mathcal{G}_{n}$, also called predictive distribution. In $P_{n}(\omega,A)$, we will often not highlight the dependence on $\omega$, and simply write $P_{n}(A)$, and write $(P_n)_{n \geq 0}$ for the sequence of random probability measures $(P_n(\cdot,\cdot))_{n \geq 0}$. Finally, for sequences of r.v.s. we denote $\mathbb{E}_{n} (X_{i}):= \mathbb{E}(X_{i}|\cG_n)$, while for   arrays $\mathbb{E}_{j} (Y_{n,i}):= \mathbb{E}(Y_{n,i}|\mathcal{F}_{n,j})$, where the r.v. and $\sigma$-field are always taken from the same row.

\section{$\xi$-a.c.i.d. Random Variables} \label{sec:acid}

A sequence of r.v.s $X=(X_{n})_{n\geq 1}$ is c.i.d., given its natural filtration $\mathcal{G}$, if for all $n\geq 1$, $X_{n+1}$ and $X_{n+2}$ are identically distributed, conditionally on $\mathcal{G}_{n}$. Equivalently, $X$ is c.i.d. if its sequence of predictive distributions $(P_{n})_{n\geq 0}$ is a measure-valued $\mathcal{G}$-martingale, i.e., $(P_{n}(A))_{n\geq 0}$ is a $\mathcal{G}$-martingale on $[0,1]$ for all $A \in \mathcal{E}$. 

Given a sequence of (possibly random $\mathcal{G}$-adapted) parameters $\xi=(\xi_{n})_{n\geq 0}$, $\xi$-a.c.i.d. r.v.s are natural generalizations of c.i.d. r.v.s, in which the sequence of predictive distributions $(P_{n}(A))_{n\geq 0}$ is an almost $\mathcal{G}$-supermartingale for all $A\in \mathcal{E}$ \citep{robbins1971convergence}, see Appendix~\ref{sec:almost.supermartingales} for a brief review. 
\begin{definition}[$\xi$-a.c.i.d. Random Variables] \label{def:acid.pred}
    Given a sequence of $\mathcal{G}$-adapted non-negative r.v.s. $\xi=(\xi_{n})_{n\geq 0}$, a sequence of r.v.s. $(X_{n})_{n\geq 1}$ is said to be $\xi$-\textit{Almost Conditionally Identically Distributed ($\xi$-a.c.i.d.)} 
     if the sequence
  $  \left( P_{n}(A) \right)_{n\geq 0}$
    forms an almost $\mathcal{G}$-supermartingale  for all $A\in \mathcal{E}$, i.e., for all $n \geq 0$ and $A\in \mathcal{E}$, a.s.
\begin{align} \label{eq:acid.pred}
    \mathbb{E}_{n}\left( P_{n+1}(A) \right) & \leq P_{n}(A) +\xi_{n}.
\end{align}
\end{definition}
Under the assumption that $\mathcal{E}$ is countably generated (e.g., when $(E,\mathcal{E})$ is Polish), an equivalent definition of $\xi$-a.c.i.d. sequences can be given in terms of TV bound between $\mathbb{P}(X_{n+1}\in \cdot |\mathcal{G}_{n})$ and  $\mathbb{P}(X_{n+2}\in \cdot |\mathcal{G}_{n})$. Specifically,

\begin{definition}[A.c.i.d., Law] \label{def:acid.law}
Given a sequence of $\mathcal{G}$-adapted non-negative r.v.s. $\xi=(\xi_{n})_{n\geq 0}$, a sequence of r.v.s. $(X_{n})_{n\geq 1}$ is said to be $\xi$-a.c.i.d. 
if for all $n\geq 0$, a.s.
\begin{equation} \label{eq:TV.acid}
    \text{TV}\left( \mathbb{P}(X_{n+2}\in \cdot | \mathcal{G}_{n}), \mathbb{P}(X_{n+1}\in \cdot | \mathcal{G}_{n})\right) \leq \xi_{n}.
\end{equation}
\end{definition}

The proof that the two definitions are equivalent is in the Appendix, subsection~\ref{sec:auxiliary}.

A trivial remark is that if $X$ is $\xi$-a.c.i.d., then it is also $\xi'$-a.c.i.d., for every sequence  $\xi'$ such that $\xi_{n}\leq \xi'_{n}$ for all $n\geq 1$. Also, if $X$ is $\xi$-a.c.i.d. for $\xi$ such that $\xi_{n}=0$ for all $n\geq 0$, then $X$ is c.i.d.

$\xi$-a.c.i.d. r.v.s were studied in \cite{bc25}, where it was shown that the integrability of $\sum_{n\geq 0} \xi_n$ implies the asymptotic exchangeability of $X$. In Theorem~\ref{thm:asym.exch}, we establish that the weaker condition $\sum_{n\geq0} \xi_n<\infty$ a.s. is sufficient to guarantee asymptotic exchangeability. In particular, when $\xi$ is deterministic, it suffices that the corresponding sequence is summable.
Some statistical examples that are  $\xi$-a.c.i.d. and, for a suitable $\xi$, asymptotically exchangeable are:

\begin{example}[Exchangeable r.v.s with noise]\label{ex:exch.with.noise} 
    
Consider a sequence $X$ defined for $n\geq 1$ as
\begin{align*}
  X_n= R_n + \epsilon_n ,\quad   R_{n}\, |\, \mu  \stackrel{iid}{\sim} \N(\mu,\sigma^{2}_{r}), \quad \epsilon_{n}  \stackrel{ind}{\sim} \N(0,\sigma^{2}_{n}), \quad \text{and} 
  \quad\mu \sim \N(0,\sigma^2_\mu)
  ,
\end{align*}
where $(\sigma^{2}_{r},\sigma^{2}_{\mu})$ are known positive constants and $(\sigma_{n}^{2})_{n\geq 1}$ is known, deterministic, and nonconstant, and $(R_n)_{n\geq 1}$ is independent of $(\epsilon_n)_{n\geq 1}$. We can interpret this as a Bayesian model in which a statistician aims to infer the unknown mean $\mu$ under heteroskedastic observation noise. The distinction between $R_n$ and $\epsilon_n$ allows us to separate the intrinsic variability of the underlying exchangeable sequence $(R_n)_{n\geq 1}$ from the heteroskedasticity induced by the noise. 

The observed sequence $X=(X_{n})_{n\geq 1}$ is neither exchangeable nor c.i.d. whenever $(\sigma^2_n)_{n\geq 1}$ is noncostant. However, it is $\xi$-a.c.i.d. whenever $(\sigma_{n}^{2})_{n\geq 1}$ is a bounded variation sequence.
Indeed, the conditional  distribution of $\mu$ given $\mathcal{G}_n$ is $\text{N}(\mu_{n},\gamma_{n})$, where $\gamma_{n}=\left(\frac{1}{\sigma^{2}_{\mu}} + \sum_{i=1}^{n}\frac{1}{\sigma_{i}^{2}+\sigma_{r}^{2}}\right)^{-1}$ and $\mu_{n} =\gamma_{n}\left(\sum_{i=1}^{n}\frac{X_{i}}{\sigma_{i}^{2}+\sigma_{r}^{2}}\right)$. Then, using the total law of probability,  $X_{n+1}|\mathcal{G}_{n} \sim \text{N}(\mu_{n}, \gamma_{n}+\sigma_{n+1}^{2}+\sigma_{r}^{2})$ and $X_{n+2}|\mathcal{G}_{n} \sim \text{N}(\mu_{n}, \gamma_{n}+\sigma_{n+2}^{2}+\sigma_{r}^{2})$. Using Theorem 1.3 of \cite{devroye2018total} to upper bound the TV distance between Gaussian distributions, we have that
\begin{equation*}
     \text{TV}\left( \mathbb{P}(X_{n+2}\in \cdot | \mathcal{G}_{n}), \mathbb{P}(X_{n+1}\in \cdot | \mathcal{G}_{n})\right) \leq \frac{3}{2}\frac{|\sigma_{n+1}^{2}-\sigma_{n+2}^{2}|}{\sigma_{n+1}^{2}+\sigma_{r}^{2}+\gamma_{n}} =: \xi_{n}.  
\end{equation*}
Many natural choices of $(\sigma^{2}_{n})_{n\geq 1}$ lead to summable $(\xi_{n})_{n\geq 0}$. In particular, since $|\sigma_{n+1}^{2}-\sigma_{n+2}^{2}|/(\sigma_{n+1}^{2}+\sigma_{r}^{2}+\gamma_{n}) \leq |\sigma_{n+1}^{2}-\sigma_{n+2}^{2}|/\sigma_{r}^{2}$, summability of $(\xi_{n})_{n\geq 0}$ holds whenever $(\sigma^{2}_{n})_{n\geq 1}$ is a bounded variation sequence (e.g., any monotone convergent sequence). Under this assumption, by Theorem~\ref{thm:asym.exch}, $X$ is asymptotically exchangeable.
\end{example}

\begin{example}[Kernel Polya Sequence]\label{ex:kernel} 

In this example, we consider recursive kernel estimators, a standard class of estimators in statistics (see, e.g., \cite{silverman2018density}). Unlike their standard use in density estimation, we use them to define a sequence of predictive distributions $(P_n)_{n\geq 1}$. There is also a substantial literature on recursively defined predictive sequences based on kernels (see, among others, \cite{mailler2017measure,berti2023kernel,sariev2024characterization,lpr26}). An important distinction from this literature is that we use the term \textit{kernel} to refer to a function $k:\mathbb{R}\mapsto [0,\infty)$ satisfying $\int_{\mathbb{R}} k(x)\di x=1$, rather than to a transition kernel. Accordingly, we follow the terminology commonly used in statistics for kernel density estimation, rather than the more general terminology used in probability theory. We return to this distinction below.
Let $(h_n)_{n\geq 0}$ be a positive sequence, called the \textit{bandwidth}, with $h_0=1$ and that depends only on $n$. Denote $k_n(x):= \frac{1}{h_n}k(x/h_n)$, and, for a set $A$ and some random variable $X$, let $K_n(A):=\int_A k_n(x) \di x$ be the corresponding measure, and $K_n(A;X):=\int_A k_n(x-X) \di x$ a its recentered version. \\
Assume that the marginal distribution of $X_1$ is the base measure $K_{0}$. Then, for all $n\geq 1$, the predictive distribution is given by
\begin{align} 
P_{n}(\cdot) & = \frac{\theta}{\theta+n} K_0(\cdot)+\frac{1}{\theta+n} \sum_{i=1}^nK_i (\cdot ; X_i) \nonumber \\
\label{eq:ex_kernel_recursion}
     &= P_{n-1}(\cdot) + \frac{1}{\theta+n}\left(K_{n}(\cdot ;X_{n}) - P_{n-1}(\cdot)\right).
\end{align}
The sequence generated is neither exchangeable nor c.i.d. \citep{west1991kernel,bc25,berti2023kernel,sariev2024characterization}. \cite{bc25} show that $X$ is $\xi$-a.c.i.d. for common choices of the kernel (e.g., Gaussian, uniform, triangular) and under suitable assumptions on the bandwidth. 

For ease of presentation, in the following, we focus on a case where $X$ is $\xi$-a.c.i.d.: we set $\theta=1$, $h_n=e^{-n}$ for all $n$, and $k(x)=\phi(x)$, where $\phi(x)$ denotes a standard Gaussian density.
Under this setting, Proposition 4 and Example 4 of \cite{bc25} show that, if $\xi$ is defined as $\xi_{n}:=  \frac{3 h^2_{n+1}}{2 n(n+1)} \sum_{i=1}^{n} \frac{1}{h^2_{i}}$ for all $n\geq 1$, then the sequence $(X_{n})_{n \geq 1}$ generated by~\eqref{eq:ex_kernel_recursion} is $\xi$-a.c.i.d with respect to $\cG$. Plugging $h_{n}=e^{-n}$ and using the geometric sum, $\sum_{i=1}^{n}e^{2i} = \frac{e^{2}(e^{2n}-1)}{e^{2}-1}$, we have that $\xi_{n} = \frac{3(1-e^{-2n})}{2n(n+1)(e^{2}-1)}$, which is summable. The choice $h_n=e^{-n}$ might appear restrictive, as it violates the usual condition $nh_n\to\infty$ for consistency of kernel density estimators (see, e.g., \cite{silverman2018density}). However, our objective here is not density estimation. Furthermore, subsequent works by \cite{christensen2026posterior,hilbert2026predictive} study the sequence in \eqref{eq:ex_kernel_recursion} under more general bandwidth choices.
\end{example}

The sequence $(P_n)_{n \geq 1}$ can be seen as a possible generalization of the classic Polya sequence \citep{blackwell1973ferguson}. Similar generalizations have been considered in several recent works \citep{mailler2017measure,bandyopadhyay2022new,berti2023kernel,sariev2024characterization,sariev2025sufficientness,chorbadzhiyska2025some}. Most of these works consider more general continuous updating schemes, but address different questions. The works most closely related to ours are \citep{berti2023kernel,sariev2024characterization,sariev2025sufficientness,chorbadzhiyska2025some,lpr26}. However, these works focus on transition kernels that preserve exchangeability, whereas our construction uses kernel estimators (in the sense defined above) to generate predictive distributions that need not preserve exchangeability. In this sense, it is closer to \citep{lpr26}.

\section{SLLN and asymptotic exchangeability}\label{sec:slln.ae}

We start by proving a SLLN for $\xi$-a.c.i.d. r.v.s.
\begin{theorem}\label{thm:slln} 
Let $(X_{n})_{n\geq 1}$ be a $\xi$-a.c.i.d. sequence, and $f:E\to \mathbb{R}$ a bounded measurable function. If $\sum_{n=0}^\infty  \xi_n  < \infty$ a.s., then there is a bounded random variable $V_{f}$ such that
\begin{equation}
    \frac{1}{n}\sum_{i=1}^{n}f(X_{i}) \to V_{f}, 
\end{equation}
a.s. and in $L^1$.
\end{theorem}
\begin{proof}
First, we show that $\int_{E}f(x)P_{n}(\di x) \stackrel{a.s.}{\to} V_{f} 
$, for some random variable $V_{f}$. First, note that $(\int_{E}f(x)P_{n}(\di x))_{n \geq 0}$ 
is an almost supermartingale. To see this, 
\begin{align}\label{eq:SLLN.formula2}
   \Bigl| & \mathbb{E}\left[ \int_{E}  f(x) P_{n+1}(\di x) | \mathcal{G}_{n}\right] - \int_{E}  f(x) P_{n}(\di x) \Bigr| \nonumber \\
   & = \left| \int_{E} f(x) ( \mathbb{P}(X_{n+2}\in \di x |\mathcal{G}_{n}) -  \mathbb{P}(X_{n+1}\in \di x |\mathcal{G}_{n})) \right| \nonumber \\
   & \leq 2 \lVert f \rVert_{\infty} TV( \mathbb{P}(X_{n+2}\in \cdot|\mathcal{G}_{n}), \mathbb{P}(X_{n+1}\in \cdot |\mathcal{G}_{n})) \nonumber \\
   & \leq 2 \lVert f \rVert_{\infty} \xi_{n},
\end{align}
where the equality follows by the tower property, the first inequality follows from the TV upper bound $|\int f Q - \int f Q'|\leq 2 \lVert f\rVert_{\infty} TV(Q,Q')$, valid for any two probability measures $Q,Q'$ and bounded function $f$, with $\lVert f\rVert_{\infty}$ denoting the uniform norm, and the last inequality follows from the $\xi$-a.c.i.d. assumption, Definition~\ref{def:acid.law}. The sequence $2 \lVert f \rVert_{\infty} \xi_{n}$ is a.s. summable. It follows that $(\int_{E}f(x)P_{n}(\di x))_{n \geq 0}$ converges a.s. by Theorem 1 of \cite{robbins1971convergence} (see Theorem~\ref{th:robsig}). 


Next, let us define
   \begin{equation*}
       Z_{n}=\sum_{i=0}^{n-1}\frac{f(X_{i+1})-\int_{E}f(x) P_{i}(\di x) }{i+1}.
   \end{equation*}
   It is easy to check that $(Z_{n})_{n\geq 0}$ is a martingale bounded in $L^{2}$, hence $Z_{n}$ converges a.s., and, by Kronecker's lemma, 
   \begin{equation} \label{eq:conv1}
       \frac{1}{n} \sum_{i=0}^{n-1}(f(X_{i+1})-\int_{E}f(x)P_{i}(\di x) ) \stackrel{a.s.} {\to} 0. 
   \end{equation}
   From $\int_{E}f(x)P_{n}(\di x) \stackrel{a.s.} {\to} 
   V_{f}$, $\frac{1}{n}\sum_{i=0}^{n-1}\int_{E}f(x)P_{i}(\di x)  \stackrel{a.s.} {\to} 
   V_{f}$. Combining it with \eqref{eq:conv1}, $ \frac{1}{n}\sum_{i=1}^{n}f(X_{i}) \stackrel{a.s.} {\to} 
   V_{f}$ and also in $L^1$ by bounded convergence.
\end{proof}

A remark is that, for c.i.d. r.v.s, Theorem 2.2 of \cite{berti2004limit} proves the SLLN for all $f \in \mathcal{H}:=\{ f:E\to \mathbb{R} \ | \ \mathbb{E}(|f(X_{1})|)<\infty\}$. The following counter-example shows that the condition $f\in\mathcal{H}$ is not enough for a SLLN for $\xi$-a.c.i.d. r.v.s and additional assumptions, such as boundedness of $f$, are needed. 
\begin{example}
Let $(X_n)_{n\ge1}$ be independent r.v.s on $\mathbb{R}$, with $X_{1}=0$, and $X_{n}\sim (1-\frac{1}{n})\delta_{0}(\cdot) +\frac{1}{n} G(\cdot)$ for $n\geq 2$, where $G$ is the distribution with density $g(y)=(y \log^2 y)^{-1} \mathbb{I}(y>e)$, $e$ being Euler's number. The sequence $(X_{n})_{n\geq 1}$ is $\xi$-a.c.i.d. with 
\begin{align}\label{eq:example.fnotbounded}
    \text{TV}(  \mathbb{P}(X_{n+1}\in \cdot|\mathcal{G}_{n}),\mathbb{P}(X_{n+2}\in \cdot |\mathcal{G}_{n})) &= \text{TV}(\mathcal{L}(X_{n+1}),\mathcal{L}(X_{n+2})) \nonumber \\
   & =\left|\frac{1}{n+1}-\frac{1}{n+2}\right|= \frac{1}{(n+1)(n+2)}=: \xi_{n},
\end{align}
where the second equality follows because $\delta_{0}$ and $G$ are mutually singular. The sequence $\xi$, defined in \eqref{eq:example.fnotbounded}, is summable. Let us consider the function $f(x)=x$, which is in $\mathcal{H}$, because $\mathbb{E}(|f(X_{1})|)=0$. Fix an integer $K\geq 1$. For $n\geq 3$, we have $Kn> e$, and therefore $$\mathbb{P}(X_{n}>Kn) =\frac{1}{n}\int_{Kn}^{\infty} g(x)dx  =\frac{1}{n}\int_{Kn}^{\infty} \frac{1}{x \log^{2}x } dx =\frac{1}{n}\int_{\log(Kn)}^{\infty} \frac{1}{y^{2}} dy = \frac{1}{n\log(Kn)}.$$

Since $\log(Kn)\leq 2\log n$ for $n\ge K$, we get
$\sum_{n\geq 3}\mathbb{P}(X_n>Kn)=\infty$. The $X_n$ are independent, so by the second
Borel--Cantelli lemma $X_n/n>K$ infinitely often a.s. Intersecting over $K\in\mathbb{N}$,
$\limsup_n X_n/n=\infty$ a.s. Finally, since $X_n\ge0$, $\frac{1}{n}\sum_{i\leq n}X_i\ge X_n/n$,
hence $\limsup_n\frac{1}{n}\sum_{i\leq n}X_i=\infty$ a.s., and the SLLN fails for $f$.

This example also shows that the condition  $\int_{E}|f(x)|\tilde{F}(dx)<\infty$ a.s. is not sufficient to guarantee the SLLN for $f$. Indeed, the directing measure $\tilde{F}$ is deterministic and equal to $\delta_{0}$. Hence, taking $f(x)=x$ as before, we have that $\int_{E}|f(x)|\tilde{F}(dx)=\int_{\mathbb{R}}|f(x)|\delta_{0}(dx)=0<\infty$. Nevertheless, as shown above, the SLLN fails for this choice of $f$.
\end{example}

Theorem 1 of \cite{bc25} establishes that $\xi$-a.c.i.d. r.v.s are asymptotically exchangeable under the condition $\mathbb{E}(\sum_{n\geq 0} \xi_{n})<\infty$. The proof relies on a generalization of Lemma 7.14 of \cite{ald85}, together with the maximal inequality for almost supermartingales (Proposition 2 of \cite{robbins1971convergence}). The following theorem strengthens this result by replacing the condition $\mathbb{E}(\sum_{n\geq 0} \xi_{n})<\infty$ with the weaker assumption $\sum_{n\geq 0} \xi_{n}<\infty$ a.s. Moreover, it establishes stable convergence of the entire shifted process, which in turn implies asymptotic exchangeability. Its statement and proof generalize Theorem 2.5 of \cite{berti2004limit} to $\xi$-a.c.i.d. sequences. The final conclusion of Theorem~\ref{thm:asym.exch}, namely asymptotic exchangeability, could be also established without (ii), for example by appealing to results in \cite{ald85,lpr26}.


\begin{theorem}\label{thm:asym.exch}
Let $(X_n)_{n\geq1}$ be a
$\xi$-a.c.i.d.\ sequence with $\sum_{n\geq0}\xi_n<\infty$ a.s. Then there exists a random probability
measure $\tilde{F}$ on $E$ such that:
\begin{enumerate}
  \item[(i)] $P_n\to\tilde{F}$ weakly, a.s.;
  \item[(ii)] $(X_n,X_{n+1},\dots)$ converges stably with representing measure $\tilde{F}^\infty$.
\end{enumerate}
In particular, $(X_n,X_{n+1},\dots)$ converges in distribution to the exchangeable law
$\mu(\cdot)=\mathbb{E}[\tilde{F}^\infty(\cdot)]$, that is, $(X_n)_{n\geq1}$ is asymptotically exchangeable.
\end{theorem}

\begin{proof}
(i) By the first part of the proof of Theorem~\ref{thm:slln}, $\int_E f \di P_n$ converges a.s. for every
bounded measurable $f$. Since $\int_{E} f \di P_n=\mathbb{E}[f(X_{n+1})\mid\mathcal{G}_n]$ and
$\sigma(X_1,X_2,\dots)=\mathcal{G}_\infty$, Lemma 2.4 of \cite{berti2004limit} (recalled as Lemma~\ref{lem:BPR.2.4}), applied with $S=E$ and $Z_n=X_n$,
yields a $\mathcal{G}_\infty$-measurable random probability measure $\tilde{F}$ on $E$ such that, for every $f\in C_b(E)$,
\begin{equation} \label{eq:ae.marginal}
  V_f=\int_E f\,d\tilde{F}\qquad\text{a.s.} 
\end{equation}
with $V_{f}$ as in Lemma~\ref{lem:ae.prod}. Letting $f$ range over a countable convergence-determining class gives $P_{n}\to \tilde{F}$ weakly a.s. 

(ii) For $n\geq1$, let $Z_n=(X_n,X_{n+1},\dots)$, an $E^\infty$-valued random variable. We aim to show that the sequence $(Z_{n})_{n\geq 1}$ converges stably with representing measure $\tilde{F}^{\infty}$, and, in particular, converges in distribution to the exchangeable law $\mu$. For this purpose, we aim at proving convergence of the finite-dimensional laws, by applying the functional monotone class theorem, Corollary~\ref{cor:fmct}.
Let $D_0\subset C_b(E)$ be a countable convergence-determining class containing the constant $1$
and closed under multiplication (e.g. the multiplicative closure of a countable
convergence-determining class). Let $\nu_n:\Omega\times\mathcal{B}(E^\infty)\to[0,1]$ be a regular
version of the conditional distribution of $(X_{n+1},X_{n+2},\dots)$ given $\mathcal{G}_n$, and let
$\nu^{(k)}_n(\omega)$ denote the projection of $\nu_n(\omega,\cdot)$ on the first $k$ coordinates,
which is a regular version of the conditional distribution of $(X_{n+1},\dots,X_{n+k})$ given
$\mathcal{G}_n$. Since $D_0$ is countable, by the definition of $\nu^{(k)}_n$, Lemma~\ref{lem:ae.prod} and
\eqref{eq:ae.marginal}, there is $\Omega_0\in\mathcal{A}$ with $\mathbb{P}(\Omega_0)=1$ such that,
for every $\omega\in\Omega_0$, all $n,k\geq1$ and all $f_1,\dots,f_k\in D_0$,
\begin{align}\label{eq:ae.fdd}
\int_{E^k}\prod_{j=1}^k f_j(x_j)\,\nu^{(k)}_n(\omega)(\di x)
 & =\mathbb{E}\Big[\prod_{j=1}^k f_j(X_{n+j})\,\Big|\,\mathcal{G}_n\Big](\omega) \nonumber \\
&  \to\prod_{j=1}^k V_{f_j}(\omega)=\prod_{j=1}^k\int_E f_j\,\di \tilde{F}(\omega).
\end{align}
Fix $\omega\in\Omega_0$ and $k\geq1$. Taking $f_i\equiv1$ for $i\neq j$ in \eqref{eq:ae.fdd}, the
$j$-th marginal of $\nu^{(k)}_n(\omega)$ converges weakly to $\tilde{F}(\omega)$. Weakly convergent
sequences of probability measures on a Polish space are tight \citep{Bil99}, so each marginal
sequence is tight on $E$ and, therefore, so is $(\nu^{(k)}_n(\omega))_{n\geq1}$ on $E^k$. Let $Q$ be
any weak limit point of $(\nu^{(k)}_n(\omega))_{n\geq1}$. Since the functions
$x\mapsto\prod_{j=1}^k f_j(x_j)$ are bounded and continuous, \eqref{eq:ae.fdd} gives, for all  $f_1,\dots,f_k\in D_0$,
\begin{equation*}
\int_{E^k}\prod_{j=1}^k f_j(x_j)\,Q(\di x)=\prod_{j=1}^k\int_E f_j\,\di \tilde{F}(\omega)
  =\int_{E^k}\prod_{j=1}^k f_j(x_j)\,\tilde{F}(\omega)^k(\di x).
\end{equation*}
By Lemma~\ref{lem:gen.class}, the class $\mathcal{C}_k:=\{x\mapsto\prod_{j=1}^k f_j(x_j): f_1,\dots,f_k\in D_0\}$ is closed under multiplication and generates
$\mathcal{B}(E^k)$, and $Q$ and $\tilde{F}(\omega)^k$ agree on it. Hence, by Corollary~\ref{cor:fmct},
$Q=\tilde{F}(\omega)^k$, and therefore $\nu^{(k)}_n(\omega)\to\tilde{F}(\omega)^k$ weakly, for all
$k\geq1$ and $\omega\in\Omega_0$.

Since on $E^\infty$ with the product topology weak convergence is equivalent to convergence of the
finite-dimensional distributions (see, e.g., \cite{Bil99}), and the finite-dimensional projections
of $\nu_n(\omega)$ converge weakly to those of $\tilde{F}(\omega)^\infty$, we have
$\nu_n(\omega)\to\tilde{F}(\omega)^\infty$ weakly for every $\omega\in\Omega_0$. That is, for every
$g\in C_b(E^\infty)$,
\begin{equation}\label{eq:ae.cond}
\mathbb{E}\big[g(X_{n+1},X_{n+2},\dots)\,\big|\,\mathcal{G}_n\big]\to\int_{E^\infty}g\,\di \tilde{F}^\infty
  \qquad\text{a.s.}
\end{equation}

Finally, we apply Lemma~\ref{lem:BPR.2.4} with $S=E^\infty$ and $Z_n=(X_n,X_{n+1},\dots)$. Its
assumptions hold: $\sigma(Z_1,Z_2,\dots)=\sigma(X_1,X_2,\dots)=\mathcal{G}_\infty$, and, by
\eqref{eq:ae.cond}, $\mathbb{E}[g(Z_{n+1})\mid\mathcal{G}_n]$ converges a.s.\ for every
$g\in C_b(E^\infty)$. Hence $(Z_n)_{n\geq1}$ converges stably with some representing measure $\gamma$,
and $\mathbb{E}[g(Z_{n+1})\mid\mathcal{G}_n]\to\int g\,\di \gamma$ a.s. Comparing with \eqref{eq:ae.cond} on a
countable convergence-determining class of $E^\infty$ gives $\gamma=\tilde{F}^\infty$ a.s., which proves (ii).
Taking $H=\Omega$ in Definition~\ref{def:stable} yields convergence in distribution to
$\mu=\mathbb{E}[\tilde{F}^\infty(\cdot)]$, which is exchangeable as a mixture of i.i.d. laws.
\end{proof}

\section{Absolute Continuity of $\tilde{F}$}   \label{sec:abs.cont}

A question of interest about the random measure $\tilde{F}$ is whether its realizations are a.s. absolutely continuous with respect to a given $\sigma$-finite measure $\lambda$. This question is of interest in Bayesian predictive inference, where $\tilde{F}$ can be interpreted as the directing measure of $(X_n)_{n \geq 1}$. It is therefore of interest to characterize the underlying $\tilde{F}$ implicitly induced by a sequence $(P_n)_{n\geq 0}$. This problem was studied for c.i.d. r.v.s in \cite{bpr13}. 

For ease of presentation, in this section, we assume $\xi$ to be non-random. 

\begin{lemma}\label{lem:tv_bound}
    Let $X=(X_n)_{n\geq 1}$ denote a $\xi$-a.c.i.d. sequence with $\sum_{n=0}^{\infty} \xi_{n}<\infty$, and $\tilde{F}$ from Theorem~\ref{thm:asym.exch}. Then, for all $n\in \mathbb{N}$,
    \begin{equation}
        \text{TV}(\mathbb{E}(\tilde{F}|\cG_{n}),P_{n}) \leq \sum_{i=n}^{\infty} \xi_i  \ \ \ \ \ \ \ \  \text{a.s.}
    \end{equation}
\end{lemma}
\begin{proof} Fix $n$ and consider the sequence $(\mathbb{E}(P_{n+M}|\cG_{n}))_{M \geq 0}$. For all $j'>j \geq 0$, it holds that
\begin{equation}\label{eq:lem_tv_finite_bound}
\text{TV} (\mathbb{E}(P_{n+j}|\cG_{n}),\mathbb{E}(P_{n+j'}|\cG_{n})) \leq \sum_{i=j}^{j'} \xi_{n+i} \ \ \ \ \ \ \ \  \text{a.s.}
\end{equation}
The sequence $(\mathbb{E}(P_{n+M}|\cG_{n}))_{M \geq 0}$ is Cauchy with respect to the TV distance: for all $\epsilon>0$, choose $N$ such that $\sum_{n=N}^{\infty} \xi_{n}< \epsilon$, then for all $j',j \geq N$, $\text{TV} (\mathbb{E}(P_{n+j}|\cG_{n}),\mathbb{E}(P_{n+j'}|\cG_{n})) \leq \epsilon$. Since the space of probability measures equipped with the TV distance is complete, $(\mathbb{E}(P_{n+M}|\cG_{n}))_{M \geq 0}$ converges in TV (to some limit) a.s. (event 1). From Corollary 1 of \cite{bc25}, $(\mathbb{E}(P_{n+M}|\cG_{n}))_{M \geq 0}$ converges weakly a.s. to $\mathbb{E}(\tilde{F}|\cG_{n})$ (event 2). On the intersection of these two events, $(\mathbb{E}(P_{n+M}|\cG_{n}))_{M \geq 0}$ converges in TV  to $\mathbb{E}(\tilde{F}|\cG_{n})$. To complete the argument, we take the limit as $M \to \infty$ in \eqref{eq:lem_tv_finite_bound} with $j=0$ and $j'=M$, and obtain
\begin{align*}
\text{TV} ( \mathbb{E}(\tilde{F}|\cG_{n}),P_{n}) & = \text{TV} (\lim_{M\to \infty} \mathbb{E}(P_{n+M}|\cG_{n}),P_{n}) = \lim_{M\to \infty} \text{TV} (\mathbb{E}(P_{n+M}|\cG_{n}),P_{n}) \\
& \leq \lim_{M\to \infty} \sum_{i=n}^{M} \xi_{i} = \sum_{i=n}^\infty \xi_i,
\end{align*}
where the second equality follows by continuity of the metric, and the inequality from \eqref{eq:lem_tv_finite_bound}.
\end{proof}

Theorem 1 of \cite{bpr13} provides the following necessary and sufficient conditions for the asymptotic directing measure $\tilde{F}$ of c.i.d. sequences to be a.s. absolutely continuous with respect to a $\sigma$-finite measure $\lambda$,
\begin{enumerate}
    \item[1'.] $\mathcal{L}(X_{1:n}) \ll \lambda^{n}$ a.s. for all $n\in \mathbb{N}$; 
    \item[2'.] $\text{TV}(P_{n},\tilde{F}) \to 0$ a.s.
\end{enumerate}
The following theorem extends Theorem 1 of \cite{bpr13} to the $\xi$-a.c.i.d. case. Conditions 1' and 2' are replaced by the apparently weaker Conditions 1 and 2 (Condition 1' always implies Condition 1, and 2' implies 2). We will return and provide some insights about these assumptions after Theorem~\ref{thm:abs.cont}.    

\begin{theorem}
    \label{thm:abs.cont}
    Let $X=(X_n)_{n\geq 1}$ denote a $\xi$-a.c.i.d. sequence with $\sum_{n=0}^{\infty} \xi_{n}<\infty$. Let $\lambda$ be a $\sigma$-finite measure on $(E,\mathcal{E})$. Then, $\tilde{F} \ll \lambda$ a.s. if and only if these two conditions hold
    \begin{enumerate}
    \item[1.] $\mathbb{P}\left(\{\omega \in \Omega : \forall B\in \mathcal{E} \ \text{s.t.} \ \lambda(B)=0, \lim_{n\to \infty} P_{n}(\omega)(B) = 0 \}\right)=1.$ 
    \item[2.] $\mathbb{P}\left(\{\omega \in \Omega :  \forall B\in \mathcal{E} \ \text{s.t.} \ \lambda(B)=0, \lim_{n\to \infty} P_{n}(\omega)(B)  = \tilde{F}(\omega)(B)\}\right) =1. $ 
    \end{enumerate}
\end{theorem} 
\begin{proof}
    ``If'' part: The set $A:= \{\omega \in \Omega :\forall B\in \mathcal{E} \ \text{s.t.} \ \lambda(B)=0, \lim_{n\to \infty} P_{n}(\omega)(B) = 0 \} \cap \{\omega \in \Omega : \forall B\in \mathcal{E} \ \text{s.t.} \ \lambda(B)=0, \lim_{n\to \infty} P_{n}(\omega)(B)  = \tilde{F}(\omega)(B) \}$ has probability one, under Conditions 1 and 2. For any $\omega \in A$, we have that, for every $B\in\mathcal{E}$ s.t. $\lambda(B)=0$, 
    \begin{equation*}
        \tilde{F}(\omega)(B) \stackrel{2.}{=} \lim_{n \to \infty}P_{n}(\omega)(B) \stackrel{1.}{=}0,
    \end{equation*}
    hence  $\tilde{F}(\omega) \ll \lambda$. Because $\mathbb{P}(A)=1$, the result follows. \\
    ``Only if'' part: Assume $\tilde{F}\ll \lambda$ a.s., define the density $f:\Omega \times E \to [0,\infty)$, which can be assumed $\mathcal{A}\otimes \mathcal{E}$-measurable \citep[][page 52]{DM82}, s.t. $\tilde{F}(\omega)(B)=\int_{B} f(\omega,x) \lambda(\di x)$ for all $\omega \in \Omega$ and all $B\in \mathcal{E}$.  Consider $B\in \mathcal{E}$, then 
    \begin{equation*}
        \int_{B} \mathbb{E}[f(\cdot,x)|\mathcal{G}_{n}] \lambda(\di x) = \mathbb{E}\left[ \int_{B} f(\cdot,x) \lambda(\di x) \,\Bigg| \,\mathcal{G}_{n} \right]= \mathbb{E}(\tilde{F}(B)\mid \mathcal{G}_{n})  \ \ \ \text{a.s.}
    \end{equation*}
    Taking the superior over $B\in\mathcal{E}$ and the limit as $n\to \infty$,  \begin{equation}\label{eq:tv_condexp_pn}
    \lim_{n \to \infty} \sup_{B\in\mathcal{E}}|\mathbb{E}(\tilde{F}(B)|\mathcal{G}_{n}) - P_{n}(B)| \leq  \lim_{n \to \infty} \sum_{i=n}^{\infty}\xi_{i} =0  \ \ \ \text{a.s.},
    \end{equation}
    where the upper bound holds a.s. by Lemma~\ref{lem:tv_bound}. In particular,  $\forall B \in \mathcal{E}$ s.t. $\lambda(B)=0$, it holds that $\int_{B}\mathbb{E}(f(\cdot,x) |\mathcal{G}_{n}) \lambda(\di x)$=0, and so   $\lim_{n\to\infty }\sup_{B:\lambda(B)=0}P_{n}(B)=0$ a.s., hence Condition 1 holds. \\
     Condition 2 follows directly from $\tilde{F} \ll \lambda$ and Condition 1. Indeed, let $A':=\{\omega \in \Omega : \forall B\in \mathcal{E} \ \text{s.t.} \ \lambda(B)=0, \lim_{n\to \infty} P_{n}(\omega)(B) = 0\} \cap \{\omega \in \Omega : \tilde{F}(\omega) \ll \lambda\}$ and note that it is a probability one event. Take $B\in \mathcal{E}$ s.t. $\lambda(B)=0$ and $\omega \in A'$. Then, from Condition 1, $\lim_{n \to \infty}P_{n}(B)=0$, and from  $\tilde{F} \ll \lambda$, $\tilde{F}(\omega)(B)=0$. Hence, $\lim_{n \to \infty}P_{n}(\omega)(B)=\tilde{F}(\omega)(B)$. 
\end{proof}

Under the $\xi$-a.c.i.d. assumption, it is possible to show that Condition 2
in Theorem~\ref{thm:abs.cont}  can be replaced by Condition 2', defined in terms of convergence in total variation of $P_{n}$.
\begin{corollary} \label{thm:corollary.TV}
     Let $X=(X_n)_{n\geq 1}$ denote a $\xi$-a.c.i.d. sequence with $\sum_{n=0}^{\infty} \xi_{n}<\infty$. Let $\lambda$ be a $\sigma$-finite measure on $(E,\mathcal{E})$. Then, $\tilde{F} \ll \lambda$ a.s. if and only if the following two conditions are satisfied
    \begin{enumerate}
    \item[1.] $\mathbb{P}\left(\{\omega \in \Omega : \forall B\in \mathcal{E} \ \text{s.t.} \ \lambda(B)=0, \lim_{n\to \infty} P_{n}(\omega)(B) = 0 \}\right)=1.$ 
    \item[2'.] $\text{TV}(P_n,\tilde{F}) \to 0$ a.s. 
    \end{enumerate}
\end{corollary}
\begin{proof}
The proof follows the same steps as that of Theorem~\ref{thm:abs.cont}, up to the last paragraph. To conclude and see that Condition 2' follows from $\tilde{F} \ll \lambda$ a.s. and Condition 1, note that 
\[
\text{TV}(P_{n},\tilde{F} ) \leq \text{TV}(\mathbb{E}(\tilde{F}|\mathcal{G}_{n}) , P_{n})+\text{TV}(\mathbb{E}(\tilde{F}|\mathcal{G}_{n}) , \tilde{F}) \ \ \ \ \text{a.s.}.
\]
The first term on the right-hand side converges to zero by \eqref{eq:tv_condexp_pn}. To show that the second term also vanishes a.s., note that $(\mathbb{E}(\tilde{F}\mid \mathcal{G}_{n}))_{n \ge 0}$ is a measure-valued martingale. Convergence in TV can be shown following the exact same argument of the second part of the proof of Theorem 1 in \cite{bpr13}.
\end{proof}

Under the c.i.d. assumption, Conditions 1 and 1' are equivalent.
 Indeed, Condition 1' is equivalent to the condition ``$P_{n} \ll \lambda$ a.s. for all $n$'', which in turn, under the c.i.d. assumption, is equivalent to Condition 1. This is because sequences of predictive distributions of c.i.d. r.v.s. are uniformly integrable martingales for each measurable set $B\in\mathcal{E}$, and hence for all $n\in\mathbb{N}$,
 \begin{equation} \label{eq:cond1}
     P_{n}(B)=\mathbb{E}(\tilde{F}(B)|\mathcal{G}_{n})=\mathbb{E}(\lim_{n \to \infty}P_{n}(B)|\mathcal{G}_{n}).
 \end{equation}
 Under Condition 1, if $B$ is any measurable set such that $\lambda(B)=0$, the right hand side of \eqref{eq:cond1} is a.s. equal to zero, and so is the left hand side, for all $n \in \mathbb{N}$. Hence, $P_{n}\ll \lambda$ for all $n$. 
 
 These last two remarks show that Theorem 1 of \cite{bpr13} can be restated using Conditions 1 and 2 - apparently weaker, but actually equivalent under the c.i.d. assumption - in place of Conditions 1' and 2'. However, this is not possible under the more general $\xi$-a.c.i.d. assumption, because Conditions 1 and 1' are not equivalent anymore. A simple counterexample of a $\xi$-a.c.i.d. sequence, satisfying Condition 1 but not Condition 1', is as follows. 
 \begin{example}
     Let $X=(X_{n})_{n\geq 1}$ be a sequence of independent r.v.s, $(a_{n})_{n\geq 1}$ a sequence of constants in $[0,1]$ converging to zero as $n\to\infty$, $\lambda$ the Lebesgue measure, and $G$ a fixed continuous distribution. Let  $X_{n} \sim (1-a_{n}) G(\cdot) +a_{n} \delta_{0}(\cdot)$, for all $n\in \mathbb{N}$. It is easy to show that
\begin{equation*}
    \text{TV}(\mathbb{P}(X_{n+1}\in \cdot|\mathcal{G}_{n}),\mathbb{P}(X_{n+2}\in \cdot |\mathcal{G}_{n})) = \text{TV}(\mathcal{L}(X_{n+1}),\mathcal{L}(X_{n+2})) =|a_{n+1}-a_{n+2}|=: \xi_{n}.
\end{equation*}
Choosing $(a_{n})_{n\geq 1}$ appropriately (of bounded variation, e.g. monotone) so that $\sum_{n=0}^{\infty}\xi_{n}<\infty$, the resulting sequence $X$ is $\xi$-a.c.i.d. and satisfies the assumptions of Theorem~\ref{thm:abs.cont}. It satisfies Condition 1, and $\tilde{F}=G$ is absolutely continuous with respect to $\lambda$. However, for every $n$, $P_{n} \nll \lambda$.
\end{example}
It is also straightforward to construct $\xi$-a.c.i.d. sequences satisfying Condition 1 (and 1') and such that $\tilde{F}\nll\lambda$. Hence, Condition 2 is violated.
\begin{example}
    Let $X=(X_{n})_{n\geq 1}$ be a sequence of r.v.s such that $X_{n} \, |\, \mu \sim \N(\mu,\sigma^{2}_{n})$ for all $n\in\mathbb{N}$, and $\mu\sim \N(0,\sigma_{\mu}^{2})$, where $(\sigma_{n})_{n\geq 1}$ is a monotone sequence of positive constants s.t. $\sigma^{2}_{n}\to 0$. Choosing $(\sigma_{n})_{n\geq 1}$ appropriately, $X$ is $\xi$-a.c.i.d. with summable weights, and $P_{n}\ll \lambda$ for all $n$, where $\lambda$ denotes the Lebesgue measure. However, $\tilde{F}=\delta_{\mu}(\cdot)$, which is singular continuous if $\sigma_{\mu}^{2}>0$ and discrete (a point mass at $0$) if $\sigma_{\mu}^{2}=0$. 
\end{example}

We conclude this section by checking Corollary~\ref{thm:corollary.TV} for exchangeable r.v.s with noise (Example~\ref{ex:exch.with.noise}).
\begin{example}\label{ex:exch_rvs}[Exchangeable r.v.s with noise (Example~\ref{ex:exch.with.noise} continued)]
     The predictive distribution is $X_{n+1}\,|\,\mathcal{G}_{n} \sim \text{N}(\mu_{n}, \gamma_{n}+\sigma_{n+1}^{2}+\sigma_{r}^{2})$, with $\mu_{n}$ and $\gamma_{n}$ defined in Example~\ref{ex:exch.with.noise}. It is straightforward to see that 1' holds. Let $(\sigma^{2}_{n})_{n\geq 1}$ be a sequence converging to a non-negative value $\sigma^{2}_{*}\geq 0$. By the SLLN for exchangeable r.v.s. and the Kolgomorov SLLN, it holds that $\frac{\sum_{i=1}^{n}X_{i}}{n}=\frac{\sum_{i=1}^{n}R_{i}}{n}+ \frac{\sum_{i=1}^{n}\epsilon_{i}}{n} \stackrel{a.s.}{\to}\mu$. Then, by Lemma~\ref{lem:newlemma}, it holds that $\mu_n \stackrel{a.s.}{\to}\mu$. Since $\gamma_{n}\to 0$ and $\sigma^{2}_{n}\to \sigma^{2}_{*}$, we have that that $\text{TV}(\N(\mu_{n},\gamma_{n}+\sigma_{n+1}^{2} +\sigma^{2}_{r}),\N(\mu, \sigma^{2}_{*}+ \sigma^{2}_{r}))\to 0$. In turns, this means that Condition 2' holds, and by Corollary~\ref{thm:corollary.TV}, the a.s. weak limit of $P_n$ is absolutely continuous.
\end{example}
 

\section{CLTs for $\xi$-a.c.i.d. sequences} \label{sec:CLTs}

In this section, three CLTs for $\xi$-a.c.i.d. r.v.s are derived. For a fixed $A \in \mathcal{E}$, Theorem~\ref{thm:clt_pred_limit} studies a sequence of predictive distributions corresponding to a sequence of $\xi$-a.c.i.d. r.v.s, centered at the asymptotic directing measure $\tilde{F}(A)$. CLTs for sequences of predictive distributions centered at $\tilde{F}(A)$ have attracted interest in Bayesian predictive inference \citep{fp25,bdlpr25}. We return to this point after presenting the theorem. Theorem~\ref{thm:clt_empirical_limit} provides the counterpart in which, instead of predictive distributions, the corresponding sequence of empirical measures is considered. Again, the centering is given by the limiting value $\tilde{F}(A)$. Finally, Theorem~\ref{thm:clt_empirical_predictive} presents an intermediate result, addressing the merging between the sequence of predictive distributions and the corresponding sequence of empirical measures. This has direct applications in Bayesian discrete-time filtering, where the predictive distribution may not be available in closed form. 

Throughout the section, conditions in Assumption~\ref{ass:clts} are assumed. 

\begin{assumption}\label{ass:clts} Given $A\in\mathcal{E}$, the sequence $(P_{n}(A))_{n\geq 0}$ satisfies the following two conditions:
\begin{align}
    \sum_{n=1}^{\infty} \mathbb{E}( n^2 (P_{n}(A)-P_{n-1}(A))^4) & < \infty, \label{eq:cltpredlimit_ass_moments}  \\
    \mathbb{E}_{n-1}(n^2 (P_{n}(A)-P_{n-1}(A))^2) & \to \eta_{0}^{2}, \hspace{0.5cm} \text{a.s.} \ \text{and} \ \text{L}^{1}, \label{eq:cltpredlimit_ass_conv_incre}
\end{align}
 for some integrable random variable $\eta^{2}_{0}$. 
\end{assumption}

Let us consider, for $A\in \mathcal{E}$ and $n\geq 0$, 
\begin{equation*}
B_{n,A}:=\sqrt{n}\left(  P_n(A) - \tilde{F}(A)\right).
\end{equation*}

\begin{theorem}\label{thm:clt_pred_limit} 
Let $(X_n)_{n \geq 1}$ be a $\xi$-a.c.i.d. sequence 
with $(\xi_{n})_{n\geq 0}$ satisfying 
\begin{equation}\label{eq:cltpredlimit_ass_xi}
\sqrt{n}\sum_{i=n}^\infty  \xi_{i}  \stackrel{a.s.}{\to} 0.
\end{equation}
If Assumption~\ref{ass:clts} and the following condition hold,
\begin{equation}\label{eq:cltpredlimit_ass_integrable}
    \sup_n \mathbb{E}( \sup_{i\geq n} n |P_{i-1}(A)-P_i(A)|) < \infty,
\end{equation}
 then $B_{n,A}\stackrel{d}{\to}Z$, where $Z$ has characteristic function $\mathbb{E}\left(\exp(-\frac{1}{2}\eta_{0}^{2}t^{2}) \right)$, with $\eta_{0}$ the limit r.v. in \eqref{eq:cltpredlimit_ass_conv_incre}.
\end{theorem}

\begin{proof} First, note that $\eqref{eq:cltpredlimit_ass_xi}$ implies $\sum_{n=0}^\infty  \xi_n  < \infty$ a.s. (e.g., see \cite{rudin1964principles} Theorem 3.22), so Theorem~\ref{thm:asym.exch} holds. Now, define the array $\{Y_{n,i},\cF_{n,i}; n \geq 0, i \geq 0 \}$ with, for all $n$, $Y_{n,0}=0$, $\cF_{n,0}=\cG_n$, and for $i\geq 1$,
\begin{align*}
    Y_{n,i} & = \sqrt{n}\left( P_{n+i-1}(A) - P_{n+i}(A) \right) 
\end{align*}
and $\mathcal{F}_{n,i}=\mathcal{G}_{n+i}$. For all $k\geq 1$, let 
\begin{align*}
    S_{n,k}=\sum_{i=1}^{k} Y_{n,i} = \sqrt{n}\left(P_n (A) - P_{n+k}(A) \right).
\end{align*}
From $P_{n}(A)\stackrel{\text{a.s.}}{\to}\tilde{F}(A)$ (by Theorem~\ref{th:robsig}), we have that $S_{n,k} \stackrel{a.s.}{\to} S_{n,\infty} =  \sqrt{n} ( P_n (A)  - \tilde{F}(A))= B_{n,A}$. Therefore, it suffices to show that the array $\{ S_{n,i},\cF_{n,i};n\geq 1,i\geq 1 \}$ satisfies the assumptions of Theorem~\ref{thm:HH3.6}. Then, the result follows. \\
We start by checking the $L^{2}$-boundedness condition~\eqref{eq:stableCLT1_sup},
\begin{align}
    \sup_{n,k \geq 1}\mathbb{E}\left( S_{n,k}^{2} \right) & = \sup_{n}\sup_{k} \mathbb{E}\left(n \left( P_{n}(A) - P_{n+k}(A)\right)^{2}\right) \nonumber  \\
    & \leq \sup_{n}\sup_{k} \sum_{j=1}^{k} \mathbb{E}\left(n \left( P_{n+j}(A) - P_{n+j-1}(A)\right)^{2}\right) \nonumber \\
    & = \sup_{n} \sum_{j=1}^{\infty} \mathbb{E}\left(n \left( P_{n+j}(A) - P_{n+j-1}(A)\right)^{2}\right) \nonumber \\
    & = \sup_{n}   \sum_{j=n}^{\infty} \mathbb{E}\left( n \left( P_{j+1}(A) - P_{j}(A)\right)^{2}\right) \label{eq:formula.proof.3.1}
\end{align}
From the assumption of $L^{1}$ convergence in \eqref{eq:cltpredlimit_ass_conv_incre}, there exists a $N$ such that for $n>N$, 
\[
  \sum_{j=n}^{\infty}   \mathbb{E}\left(  \left( P_{j+1}(A) - P_{j}(A)\right)^{2}\right) \leq    \sum_{j=n}^{\infty} \frac{\mathbb{E}(\eta_0^2)+\epsilon}{j^2} \leq \frac{\mathbb{E}(\eta_0^2)+\epsilon}{n-1}.
\] 
Therefore, for all $n>N$, $\sum_{j=n}^{\infty} \mathbb{E}\left( n \left( P_{j+1}(A) - P_{j}(A)\right)^{2}\right) \leq 2\left(\mathbb{E}(\eta_0^2)+\epsilon\right) <\infty$ is uniformly bounded, due to the integrability of $\eta_{0}^{2}$. This implies boundedness of \eqref{eq:formula.proof.3.1}, hence condition ~\eqref{eq:stableCLT1_sup}. \\
The nesting condition~\eqref{eq:HH3.6_nestedsigmafields} holds by construction. To prove that condition~\eqref{eq:HH3.6_conv_sum2} holds, we follow the strategy of \cite{fortini2020quasi} Theorem 5. We can write 
\begin{equation*}
    \sum_{i\geq 1}Y_{n,i}^{2}=\sum_{i=1}^{\infty}n(P_{n+i-1}(A)-P_{n+i}(A))^{2}=\sum_{i=n}^{\infty}n(P_{i}(A)-P_{i-1}(A))^{2}=n\sum_{i=n}^{\infty}\frac{Q_i}{i^2},
\end{equation*}
where $Q_{i}:=i^{2}(P_{i}(A)-P_{i-1}(A))^{2}$.
Noting that $n \sum_{i=n}^\infty\frac{1}{i^2}\to 1$, under Assumption~\ref{ass:clts}, by Lemma 4.1 in \cite{crimaldi2016fluctuation} (see Lemma~\ref{lem:crimaldi_4.1} in the Appendix), applied to $(Q_n)_{n \geq 1}$,  we obtain that $\sum_{i\geq 1}Y_{n,i}^{2}\stackrel{a.s.}{\to} \eta_{0}^{2}$. \\
In order to check condition~\eqref{eq:HH3.6_max_conv}, let us rewrite
\begin{align*}
\sup_{i\geq 1} |Y_{n,i}|& = \sup_{i\geq 1} \sqrt{n} | P_{n+i-1}(A) - P_{n+i}(A)|  \\
&= \sup_{i\geq n} \sqrt{n} | P_{i-1}(A) - P_{i}(A)| \leq \sup_{i\geq n} \sqrt{i} | P_{i-1}(A) - P_{i}(A)|.
\end{align*}
From the just proven condition~\eqref{eq:HH3.6_conv_sum2}, we have that $\sum_{i\geq 1}Y_{n,i}^{2}\stackrel{a.s.}{\to} \eta_{0}^{2}$, which implies that $n (P_{n}(A)-P_{n-1}(A))^2 \stackrel{a.s.}{\to}0$. This in turn implies the a.s. convergence to $0$ of the right side of the last equation. Hence, condition~\eqref{eq:HH3.6_max_conv} holds. 
\\
Condition~\eqref{eq:HH3.6_unif_bound} is assumed in assumption~\eqref{eq:cltpredlimit_ass_integrable}. For conditions~\eqref{eq:HH3.6_conv_exp} and \eqref{eq:HH3.6_conv_exp2}, note that 
 \begin{align*}
\sum_{i=1}^{\infty}|\mathbb{E}_{i-1}\left(Y_{n,i}\right) | 
&= \sum_{i=1}^{\infty} \sqrt{n} \, | \mathbb{E}_{n+i-1}\left(P_{n+i}(A) - P_{n+i-1}(A)\right) | \leq \sum_{i=1}^{\infty} \sqrt{n} \xi_{n+i-1} =  \sqrt{n} \sum_{i=n}^{\infty}  \xi_{i}.
\end{align*}
Assumption~\eqref{eq:cltpredlimit_ass_xi} implies that the right hand side converges a.s. to zero, hence condition \eqref{eq:HH3.6_conv_exp} follows. With an analogous argument, we have that $\sum_{i=1}^{\infty}\mathbb{E}_{i-1}\left(Y_{n,i}\right)^2 \leq n \sum_{i=n}^{\infty}  \xi^2_{i}$. Condition \eqref{eq:HH3.6_conv_exp2} holds, noting that for a sequence $(a_n)_{n\geq 0}$ such  $\sum_{i=n}^\infty a_i \to 0$ and $a_i\geq 0$ for all $i$, it holds $\sum_{i=n}^\infty a^2_i \to 0$ (by the comparison test).
\end{proof}


Assumption~\ref{ass:clts} and \eqref{eq:cltpredlimit_ass_xi} are stronger than necessary. To see this, note that the proof establishes several almost-sure convergence results, whereas Theorem~\ref{thm:HH3.6} only requires convergence in probability. In fact, Theorem~\ref{thm:clt_pred_limit} could also have been stated in terms of almost-sure conditional convergence (\citep{Cr09}, see Definition~\ref{def:ascond}). This mode of stochastic convergence is stronger than convergence in distribution and (mixing) stable convergence. We nevertheless prefer to state the theorem using the more conventional language of stable convergence.

Almost-sure conditional convergence to a Gaussian limit for a suitably rescaled version of $P_n(A)-\tilde F(A)$, conditional on $\mathcal G_n$, has attracted interest in Bayesian predictive inference; see \citep{bdlpr25,fp25} for a review. Informally, the motivation is that, although $\xi$-a.c.i.d. (and c.i.d.) sequences do not admit a representation theorem, asymptotic exchangeability implies that their tails $(X_n,X_{n+1},\ldots)$ admit such a representation asymptotically as $n\to\infty$. This provides a Bayesian interpretation of the conditional distribution of $P_n(A)-\tilde F(A)$ given $\mathcal G_n$, for which a CLT yields an asymptotic Gaussian approximation. We refer to \cite{fortini2020quasi,fortini2023prediction} for further discussion of when such an approximation is justified from a Bayesian standpoint.

Theorem~\ref{thm:clt_pred_limit} provides a basis for constructing asymptotic intervals for $\tilde{F}(A)$, which can be interpreted as credible intervals under the Bayesian interpretation discussed above. Specifically, let $V_{n,A}:=\frac{1}{n}\sum_{i=1}^{n}i^{2}(P_{i}(A)-P_{i-1}(A))^{2}$ and $z_{q}$ denote the $q\%$ quantile of a Standard Normal distribution. Combining the last part of Theorem~\ref{thm:HH3.6} (see \eqref{eq:HH3.6_conclusion}) with the almost sure limit $V_{n,A} \to \eta_{0}^{2}$ guaranteed by Assumption~\ref{ass:clts} yields the following corollary to Theorem~\ref{thm:clt_pred_limit}.
\begin{corollary} \label{thm:corollary.interval}
Under the assumptions of Theorem~\ref{thm:clt_pred_limit}, a $(1-\alpha)\%$ asymptotic credible interval for $\tilde{F}(A)$ is
    \begin{equation*}
      \Big[P_{n}(A) - z_{1-\frac{\alpha}{2}}\sqrt{V_{n,A}/n}, P_{n}(A) + z_{1-\frac{\alpha}{2}}\sqrt{V_{n,A}/n}\Big]. 
    \end{equation*}
\end{corollary}

After $(P_n)_{n \geq 0}$, it is natural to characterize the limiting behavior of the empirical measure. A standard probabilistic argument can be used to show that the empirical distribution of a $\xi$-a.c.i.d. sequence with summable weights converges weakly a.s. to the same limiting random measure $\tilde{F}$. Thus, the obvious next step is to consider a CLT for such measure. As an intermediate step, we consider the deviation between the empirical and predictive distributions. Specifically, let us denote, for $A\in\mathcal{E}$ and $n\geq 0$,
\[
C_{n,A}:=\sqrt{n} \,\left(\frac{1}{n}\sum_{i=1}^{n} \mathbb{I}(X_{i}\in A)   -P_n (A) \right).
\]
\begin{theorem}\label{thm:clt_empirical_predictive}
  Let $(X_n)_{n \geq 1}$ be a $\xi$-a.c.i.d. sequence 
  with $(\xi_{n})_{n\geq 0}$  satisfying
\begin{equation}\label{eq:clt3_ass_xi}
 \frac{1}{\sqrt{n}} \sum_{i=1}^n i \, \xi_i  \stackrel{a.s.} {\to}  0.
    \end{equation}
If Assumption~\ref{ass:clts} holds, then $C_{n,A}\stackrel{d}{\to}Z$ (stably), where $Z$ has characteristic function $\mathbb{E}\left(\exp(-\frac{1}{2}\eta_2^2 t^{2}) \right)$ for some a.s. finite real random variable $\eta_2^2$.
\end{theorem}
\begin{proof} First, note that \eqref{eq:clt3_ass_xi} implies that $\xi_{n}/n^{3/2-\delta} \stackrel{a.s}{\to} 0$, for some $\delta>0$. So $\sum_{n=0}^\infty  \xi_n  < \infty$ a.s., and Theorem~\ref{thm:asym.exch} applies. 
Let us define the triangular array $\{Y_{n,i},\cF_{n,i}; n\geq0, 1\leq i \leq n\}$, where, for $n\geq 1$, $i=1,\ldots,n$, 
\[
Y_{n,i}:=\frac{1}{\sqrt{n}}(\mathbb{I}(X_{i}\in A) -i P_i (A)+(i-1)P_{i-1}(A )),
\]
and $\cF_{n,i}=\cG_i$. Let $S_{n,i}=\sum_{j=1}^{i}Y_{n,j}$, and notice that $S_{n,n}=C_{n,A}$ for all $n\geq 1$. We now show that the triangular array $\{S_{n,i}, \mathcal{F}_{n,i};n\geq 1,i\leq n\}$ satisfies the conditions of Theorem~\ref{thm:HH3.2}, which in turn proves the theorem. \\
The nesting condition~\eqref{eq:stableCT_nestedsigmafields} holds by construction. 
To prove condition \eqref{eq:stableCLT_conv_sum2}, we rewrite the sum of the squared increments as the sum of three series
\begin{align}\label{eq:clt3_squaredterms}
  \sum_{i=1}^{n} Y_{n,i}^2=& \frac{1}{n}\sum_{i=1}^{n} (\mathbb{I}(X_{i}\in A)-P_{i-1}(A))^2 +  \frac{1}{n}\sum_{i=1}^n i^2 (P_{i}(A)-P_{i-1}(A))^2\\
  &-2 \frac{1}{n}\sum_{i=1}^{n} (\mathbb{I}(X_{i}\in A)-P_{i-1}(A))i(P_{i}(A)-P_{i-1}(A)):=S_{n,1}+S_{n,2}-2 S_{n,3}.\nonumber
\end{align}
Firstly, we prove that $S_{n,1}$ converges a.s. to $\tilde{F}(A)-\tilde{F}^2(A)$. Indeed, by Theorem~\ref{thm:slln}, $\frac{1}{n}\sum_{i=1}^{n} \mathbb{I}(X_{i}\in A)\stackrel{a.s.}{\to}\tilde{F}(A)$, and, by Theorem~\ref{th:robsig}, 
$P_n(A)\stackrel{a.s.}{\to}\tilde{F}(A)$, which implies 
$\frac{1}{n}\sum_{i=1}^{n}P^2_i(A)\stackrel{a.s.}{\to}\tilde{F}^2(A)$. For the cross-product of $S_{n,1}$, 
let us define
\[
M_n:=\sum_{i=1}^n\frac{\mathbb{I}(X_{i}\in A)P_{i-1}(A)-P^2_{i-1}(A)}{i}.
\]
$(M_n)_{n \geq 1}$ is a martingale, with $\sup_n \mathbb{E} (M_n^2)<\infty$ (the sum of the squared increments is bounded), hence converges a.s. An application of the Kronecker Lemma along with the same argument used in Theorem~\ref{thm:slln} gives $\frac{1}{n}\sum_{i=1}^{n} \mathbb{I}(X_{i}\in A)P_{i-1}(A)\stackrel{a.s.}{\to}\tilde{F}^2(A)$. 

Secondly, $S_{n,2}$ converges a.s. to $\eta_{0}$. Indeed, let $Q_n:= n^2 (P_{n-1}(A)-P_n(A))^2$, so that $S_{n,2}=\frac{1}{n} \sum_{i=1}^n Q_i$. By Assumption~\ref{ass:clts} and Lemma 4.1 a) of \cite{crimaldi2016fluctuation} with $a_n=1, b_n=n$ (see Appendix Lemma~\ref{lem:crimaldi_4.1}), $S_{n,2}$ converges a.s. to $\eta_{0}^2$.
\\
Finally, the cross product $S_{n,3}$ is a.s. Cauchy. Indeed, by Cauchy-Schwartz, $|S_{n,3}|\leq (S_{n,1}\, S_{n,2})^{1/2}$, and $(S_{n,1})_{n \geq 1}$ and $(S_{n,2})_{n \geq 1}$ are a.s. convergent series. Hence, for all $\omega \in \Omega$, $S_{n,3}$ is a convergent series (Theorem 3.22 in \cite{rudin1964principles}). Let $\eta_1^2$ denote the limiting r.v.\\
To sum up, condition~\eqref{eq:stableCLT_conv_sum2} holds, and the limit is $\eta_2^2=\eta_{0}^2+\eta^2_1+\tilde{F}(A)-\tilde{F}^2(A)$. \\
For condition \eqref{eq:stableCT_max_conv}, note that $Y_{n,i}=\sqrt{\frac{i}{n}}Y_{i,i}$. Then, condition~\eqref{eq:stableCT_max_conv} holds because $\max_{i\leq n} Y_{n,i}^{2} = \max_{i\leq n} \frac{i}{n}Y_{i,i}^{2} \stackrel{a.s.}{\to} 0$. For condition \eqref{eq:stableCLT_unif_bound}, rewrite
\[
\sup_n \mathbb{E}(\max_{i\leq n} Y_{n,i}^2) \leq \sup_n \mathbb{E}(S_{n,1}) + \sup_n \mathbb{E}(S_{n,2}) + \sup_n \mathbb{E}(S_{n,3}) <\infty,
\]
where the bound for the second and third terms hold by Assumption~\ref{ass:clts}.\\
To prove conditions \eqref{eq:stableCLT_conv_exp},
\[
\sum_{i=1}^{n}|\mathbb{E}_{i-1}(Y_{n,i})|= \frac{1}{\sqrt{n}}\sum_{i=1}^{n}i|P_{i-1}(A)-\mathbb{E}_{i-1}(P_i(A))|\leq \frac{1}{\sqrt{n}}\sum_{i=1}^{n} i \xi_i \stackrel{a.s.}{\to} 0,
\]
where the a.s. limit holds by Assumption \eqref{eq:clt3_ass_xi}. Condition  \eqref{eq:stableCLT_conv_exp} holds by Theorem 3.45 in \cite{rudin1964principles}. Finally, condition~\eqref{eq:stableCLT_conv_exp2} follows from condition~\eqref{eq:stableCLT_conv_exp} by a comparison test.
\end{proof}

 It may be of interest to determine the limiting variance $\eta^2_{2}$, appearing in Theorem~\ref{thm:clt_empirical_predictive}. Assumption~\ref{ass:clts}, specifically condition \eqref{eq:cltpredlimit_ass_conv_incre}, provides a way to compute $\eta_{0}^2$. An analogous argument yields the term $\eta_1^2$, the a.s. limit of the sum $S_{n,3}$. To compute it, note that the sequence $\left( \mathbb{E}_{n-1} \left( \big( \mathbb{I}(X_{n} \in A) - P_{n-1}(A) \big) n \big( P_{n}(A) - P_{n-1}(A) \big) \right) \right)_{n \geq 1}$ converges a.s. to $\eta_1^2$. To see this, let $Z_n:=(\mathbb{I}(X_{n}\in A)-P_{n-1}(A))n(P_{n}(A)-P_{n-1}(A))$, $a_n=1$, $b_n=n$, note that 
   \[
   \sum_{n=1}^\infty \frac{\mathbb{E}(Z_n^2)}{n^2} \leq \left[\sum_{n=1}^\infty \frac{\mathbb{E}((\mathbb{I}(X_{n}\in A)-P_{n-1}(A))^2)}{n^2}\sum_{n=1}^\infty \mathbb{E}((P_{n}(A)-P_{n-1}(A))^2)\right]^{1/2}<\infty,
   \]
   where the first series is trivially bounded, while the second one is bounded by Assumption~\ref{ass:clts}, condition~\eqref{eq:cltpredlimit_ass_moments}. Therefore, the conditions of Lemma~4.1 in \cite{crimaldi2016fluctuation} are satisfied, and the claim follows.


Finally, we consider a CLT between the empirical and the asymptotic directing measure $\tilde{F}$.  For $A\in \mathcal{E}$ and $n\geq 0$, let us define 
\begin{equation*}
W_{n,A}:=\sqrt{n}\left(  \frac{1}{n}\sum_{i=1}^{n} \mathbb{I}(X_{i}\in A)  - \tilde{F}(A)\right).
\end{equation*}

\begin{theorem}\label{thm:clt_empirical_limit} 
Let $(X_n)_{n \geq 1}$ be a $\xi$-a.c.i.d. sequence, and 
assume that the assumptions of  Theorems~\ref{thm:clt_pred_limit} and \ref{thm:clt_empirical_predictive} hold. Then, $W_n\stackrel{d}{\to}Z$, where $Z$ has characteristic function $\mathbb{E}\left(\exp(-\frac{1}{2}\eta^{2}_{3}t^{2}) \right)$, for some non-negative random variable $\eta^{2}_3$. 
\end{theorem}

\begin{proof} Let us define the array $\{Y_{n,i},\cF_{n,i};n \geq 0,i\geq 0 \}$ with, for all $n$, $Y_{n,0}=0$, $\cF_{n,0}=\cG_n$, and
\begin{align*}
    Y_{n,i}=\begin{cases}
    \frac{1}{\sqrt{n}}\, (\mathbb{I}(X_{i}\in A) -i P_i (A)+(i-1)P_{i-1}(A ))  & \text{for }1\leq i \leq n\\
    \sqrt{n} \, (P_{i-1}(A)-P_i(A))   & \text{for } i > n,
    \end{cases} 
\end{align*}
and, for $i\leq n$, $\mathcal{F}_{n,i}=\mathcal{G}_{i}$  and, for $i>n$, $\mathcal{F}_{n,i}=\mathcal{G}_{n+i}$ . For $k>n$, let 
\begin{align*}
    S_{n,k}=\sum_{i=1}^{k} Y_{n,i} = \sqrt{n}\left(\frac{1}{n} \sum_{i=1}^n \mathbb{I} (X_i \in A)  - P_{k}(A) \right).
\end{align*}
By Theorem~\ref{th:robsig}, $S_{n,k}  \stackrel{a.s.}{\to} S_{n,\infty}$, where 
\begin{align*}
    S_{n,\infty}:=  \sqrt{n} \left( \frac{1}{n} \sum_{i=1}^n \mathbb{I} (X_i \in A)   - \tilde{F}(A)\right).
\end{align*}
Because $S_{n,\infty}=W_{n,A}$, it is enough to show that the infinite array $\{ S_{n,i},\cF_{n,i};n\geq 1,i\geq 1 \}$ satisfies the assumptions of Theorem~\ref{thm:HH3.6} to prove the theorem.\\
Condition~\eqref{eq:HH3.6_nestedsigmafields} holds by construction. Condition~\eqref{eq:HH3.6_conv_sum2} holds by combining the arguments used in Theorems~\ref{thm:clt_pred_limit} and \ref{thm:clt_empirical_predictive}. To see this, 
\begin{equation*}
    \sum_{i\geq 1}Y_{n,i}^{2}=\sum_{i=1}^{n}\frac{1}{n}(\mathbb{I}(X_{i}\in A) -i P_i (A)+(i-1)P_{i-1}(A ))^2+\sum_{i=n+1}^{\infty}n(P_{i}(A)-P_{i-1}(A))^{2},
\end{equation*}
which is exactly the sum of the series analyzed in the proofs of the two theorems. Let us denote by $\eta_3^2$ the a.s. limit of the series. Conditions~\eqref{eq:HH3.6_max_conv} and \eqref{eq:HH3.6_unif_bound} can be proven using the same decomposition. \\
Finally, the conditional expectation of $Y_{n,i}$ given $\mathcal{F}_{n,i-1}$ satisfies 
\begin{align*}
   |\mathbb{E }_{i-1}( Y_{n,i})|\leq \begin{cases}
    \frac{i-1}{\sqrt{n}}\, \xi_i  & \text{for }1\leq i \leq n\\
    \sqrt{n}\, \xi_i  & \text{for } i > n,
    \end{cases} 
\end{align*}
which implies that 
\[
\sum_{i=1}^{\infty} |\mathbb{E}_{i-1}( Y_{n,i})| \leq   \frac{1} {\sqrt{n}} \sum_{i=1}^n (i-1) \, \xi_i +\sqrt{n} \sum_{i=n}^\infty \xi_n \stackrel{a.s.}{\to} 0.
\]
Hence, condition~\eqref{eq:HH3.6_conv_exp} holds. Condition~\eqref{eq:HH3.6_conv_exp2} follows directly from condition~\eqref{eq:HH3.6_conv_exp}.
\end{proof}

The asymptotic variance $\eta_3^2$ can be computed by an argument analogous to the one suggested for $\eta_2^2$. Theorem~\ref{thm:clt_empirical_limit} is the counterpart of Theorem 3.1 in \cite{berti2004limit}. The latter is based on an assumption that is typically hard to check in applications, whereas our proof only requires the assumptions of Theorems~\ref{thm:clt_pred_limit} and \ref{thm:clt_empirical_predictive}, making the argument possibly more transparent and easier to verify. We illustrate it in the following example, where we check the conditions of Theorems~\ref{thm:clt_pred_limit}-\ref{thm:clt_empirical_limit}.



\begin{example}[Kernel Polya sequence (Example~\ref{ex:kernel} continued)] Setting $h_n=e^{-n}$, $\xi$ is a deterministic sequence with $\xi_{n} = \frac{3(1-e^{-2n})}{2n(n+1)(e^{2}-1)}$. It is immediate to check that $\frac{1}{\sqrt{n}}\sum_{i=1}^{n} i \xi_{i} \to 0$, hence condition~\eqref{eq:clt3_ass_xi} holds.

\noindent  Now, choose $A\in \mathcal{E}$ such that $\tilde{F}(\partial A)=0$ a.s. For Assumption~\ref{ass:clts}, condition~\eqref{eq:cltpredlimit_ass_moments} is also straightforward. Indeed, using the recursive representation~\eqref{eq:ex_kernel_recursion},
\begin{align*}
\sum_{n=1}^\infty \mathbb{E}( n^2 (P_{n}(A)-P_{n-1}(A))^4) & \leq  \sum_{n=1}^\infty \mathbb{E}\left( \frac{1}{(n+1)^2}(K_n(A,X_n)-P_{n-1}(A))^4\right) \\
& < \sum_{n=1}^\infty \frac{1}{(n+1)^2} <\infty.
\end{align*}
To check~\eqref{eq:cltpredlimit_ass_conv_incre}, using \eqref{eq:ex_kernel_recursion}, let us rewrite
\begin{align*}
\mathbb{E}_{n-1} (n^{2} &(P_{n}(A) -  P_{n-1}(A))^{2})  =\frac{n^2}{(n+1)^2}
\mathbb{E}_{n-1} ((K_n(A,X_n)-P_{n-1}(A))^2)    \\
&=\frac{n^2}{(n+1)^2} [ \mathbb{E}_{n-1} (K_n(A,X_n)^2) +P^2_{n-1}(A) -2 P_{n-1}(A)\mathbb{E}_{n-1} (K_n(A,X_n)) ].
\end{align*}
Since $A$ is a $\tilde{F}$-continuity set, by the extended continuous mapping theorem (e.g., Theorem 1.11.1 \cite{van1996weak}) applied to $\mathbb{E}_{n-1} (K_n(A,X_n))$ and $\mathbb{E}_{n-1} (K_n(A,X_n)^2)$, it holds that 
\[
\frac{n^2}{(n+1)^2}\mathbb{E}_{n-1} (K_n(A,X_n)-P_{n-1}(A))^2) \stackrel{a.s.}{\to} \tilde{F}(A) - \tilde{F}^2(A)=:\eta^{2}_{0}
\]
Being $(K_n(A,X_n)-P_{n-1}(A))^2)_{n \geq 1}$ bounded, this also implies $L^{1}$ convergence. \\
Finally, the last condition to check is \eqref{eq:cltpredlimit_ass_integrable}, which is also trivial. Indeed,
\[
    \sup_n \mathbb{E}( \sup_{i\geq n} n |P_{i-1}(A)-P_i(A)|) \leq  \sup_n \mathbb{E}( \sup_{i\geq n} |K_i(A,X_i)-P_{i-1}(A)|) \leq 1.
    \]
Now, we turn to computing the asymptotic variances in closed form. $\eta_0^2$ is the one of Theorem~\ref{thm:clt_pred_limit}.
To compute the asymptotic variances in Theorems~\ref{thm:clt_empirical_predictive}–\ref{thm:clt_empirical_limit}, one can follow the discussion after Theorem~\ref{thm:clt_empirical_predictive}. In particular, $\eta_1^2$ is obtained by considering the following limit
\begin{align*}
\mathbb{E}_{n-1}&((\mathbb{I}(X_{n}\in A)-P_{n-1}(A))n(P_{n}(A)-P_{n-1}(A))\\
& = \mathbb{E}_{n-1}((\mathbb{I}(X_{n}\in A)-P_{n-1}(A))\frac{n}{n+1}(K_n(A,X_n)-P_{n-1}(A))) \stackrel{a.s.}{\to} 0,
\end{align*}
which implies that $\eta_{2}^{2}=2(\tilde{F}(A) - \tilde{F}^2(A))$. $\eta_3^2$ can be computed via a similar argument.
\end{example}

\subsection{Multivariate CLTs} The results obtained so far can be extended to the multivariate setting, which has both applied and theoretical consequences. On the applied side, it allows one to construct credible intervals for multiple sets jointly; under the c.i.d. assumption, see \cite{fortini2020quasi} and \cite{fp25}. On the theoretical side, multivariate CLTs, such as Theorem~\ref{thm:clt_pred_limit_mult}, together with a condition on the tightness of the sequence,  allow one to extend Theorem~\ref{thm:clt_pred_limit_mult} to a Uniform CLT over the Donsker-class $\mathcal{D}=\{\mathbb{I}_{(-\infty,t]}:t\in \mathbb{R}\}$ (see Theorems 1.5.4 and 1.5.6 of \cite{van1996weak}). However, proving tightness on specific examples is usually not straightforward, and therefore the study of Uniform CLTs for $\xi$-a.c.i.d. sequences is left to future research.

 We restrict sets $A$ to be semi-open intervals, and for every $n\geq 0$ and $t\in \mathbb{R}$, let us define
\[
\Delta_{n,t}:=P_n((-\infty,t])-P_{n-1}((-\infty,t]).
\]
In a similar way, for $k\in \mathbb{N}$ and $\mathbf{t}=(t_1,\ldots,t_k)^T$, we can define the vector $\mathbf{\Delta}_{n,\mathbf{t}}=(\Delta_{n,t_1},\ldots,\Delta_{n,t_k})^T$.
\begin{assumption}\label{ass:clts_mult} For all 
$i,j=1,\ldots,k$, 
we assume the following conditions hold:
\begin{equation}\label{eq:multi_cltpredlimit_ass_moments}
    \sum_{n=1}^{\infty} \mathbb{E}( n^2 \Delta_{n,t_j}^4) < \infty,
\end{equation}
\begin{equation}\label{eq:multi_cltpredlimit_ass_conv_incre}
    \mathbb{E}_{n-1}(n^2 \Delta_{n,t_j}\Delta_{n,t_i}) \to H[j,i] \hspace{0.5cm} \text{a.s.} \ \text{and} \ \text{L}^{1},
\end{equation}
 for some positive definite matrix $H$. 
\end{assumption}

We consider a vector version of Theorem~\ref{thm:clt_pred_limit}, and also include in the statement an estimator of the asymptotic covariance matrix. 
Consider
\begin{equation*}
\mathbf{B}_{n,\mathbf{t}}:=\sqrt{n} V_{n,\mathbf{t}}^{-1/2}\left(  P_n((-\infty,t_1]) - \tilde{F}((-\infty,t_1]),\ldots, P_n((-\infty,t_k]) - \tilde{F}((-\infty,t_k])\right)^T, 
\end{equation*}
where
\[
V_{n,\mathbf{t}} = \frac{1}{n} \sum_{i=1}^n i^2 \mathbf{\Delta}_{i,\mathbf{t}}\mathbf{\Delta}_{i,\mathbf{t}}^T
\]
 is an estimator for $H$. 

\begin{theorem}\label{thm:clt_pred_limit_mult} 
Let $(X_n)_{n \geq 1}$ be a $\xi$-a.c.i.d. sequence 
with $(\xi_{n})_{n\geq 0}$  satisfying 
\eqref{eq:cltpredlimit_ass_xi}. If Assumption~\ref{ass:clts_mult} and the following condition hold
\begin{equation}\label{eq:multi_cltpredlimit_ass_integrable}
    \sup_n \mathbb{E}( \sup_{i\geq n} n |\Delta_{i,t_j}|) < \infty,
\end{equation}
for all $j=1,\ldots,k$, then $\mathbf{B}_{n,\mathbf{t}}\stackrel{d}{\to}Z$, where $Z$ has characteristic function $\mathbb{E}\left(\exp(-\frac{1}{2}\mathbf{s}^T \mathbf{s}) \right)$.
\end{theorem}
The proof proceeds along the same lines as that of Theorem~\ref{thm:clt_pred_limit}, combined with the Cramér–Wold device. This yields convergence in distribution of $B_{n,\mathbf{t}}$ to a multivariate Gaussian with covariance matrix $H$. It is then sufficient to establish the almost sure convergence of $V_{n,\mathbf{t}}$ to $H$, which follows from Proposition~2.5 of \cite{fp25}.

\begin{appendix} \label{sec:appendix}
\section{Auxiliary results}\label{sec:auxiliary}

Let $S$ be a Polish space, $\mathcal{P}(S)$ the set of probability measures on $\mathcal{B}(S)$, and
$\Sigma$ the $\sigma$-field on $\mathcal{P}(S)$ generated by the maps $p\mapsto p(B)$, $B\in\mathcal{B}(S)$.
A \emph{random measure} on $S$ is a measurable map $\gamma:(\Omega,\mathcal{A})\to(\mathcal{P}(S),\Sigma)$.

\begin{definition}[Stable convergence]\label{def:stable}
Let $(Z_n)_{n\geq1}$ be $S$-valued random variables on $(\Omega,\mathcal{A},\mathbb{P})$. We say that
$(Z_n)_{n\geq1}$ \emph{converges stably} if, for every $H\in\mathcal{A}$ with $\mathbb{P}(H)>0$, $(Z_n)_{n\geq1}$
converges in distribution under $\mathbb{P}(\cdot\mid H)$ to some probability law $\mu_H$ on $S$.
In this case, there is a random measure $\gamma$ on $S$ such that
\begin{equation*}
\mu_H(\cdot)=\int_\Omega\gamma(\omega)(\cdot)\,\mathbb{P}(d\omega\mid H)
  \qquad\text{for every } H\in\mathcal{A} \text{ with } \mathbb{P}(H)>0,
\end{equation*}  
and $(Z_n)_{n\geq1}$ is said to \emph{converge stably with representing measure} $\gamma$.
\end{definition}

\begin{definition}[Almost sure conditional convergence, \citep{Cr09}]\label{def:ascond} 
Let $(Z_n)_{n\geq1}$ be $S$-valued random variables on $(\Omega,\mathcal{A},\mathbb{P})$ and $\mathcal G$ a conditioning system. Let us denote by $K_n$ a version of the conditional distribution of $Z_n$ given $\mathcal G_n$. If $K$ is a transition kernel such that, for almost every $\omega\in\Omega$, the sequence $(K_n(\omega,\cdot))_{n \geq 1}$ of probability measures on $S$ converges weakly to the probability measure $K(\omega,\cdot)$, then we shall say that, with respect to the conditioning system $\mathcal G$, the sequence 
$(Z_n)_{n\geq1}$ converges to $K$ in the sense of  \emph{almost sure conditional convergence}. 
\end{definition}

\begin{lemma}[\cite{berti2004limit}, Lemma~2.4]\label{lem:BPR.2.4}
Let $S$ be a Polish space and $(Z_n)_{n\geq1}$ a sequence of $S$-valued random variables on
$(\Omega,\mathcal{A},\mathbb{P})$. If $\sigma(Z_1,Z_2,\dots)\subset\sigma\big(\bigcup_n\mathcal{G}_n\big)$ and
$\big(\mathbb{E}[f(Z_{n+1})\mid\mathcal{G}_n]\big)_{n\geq1}$ converges a.s.\ for each $f\in C_b(S)$, then there
is a random measure $\gamma$ on $S$ such that $(Z_n)_{n\geq1}$ converges stably with representing
measure $\gamma$ and
\begin{equation*}
  \mathbb{E}[f(Z_{n+1})\mid\mathcal{G}_n]\to\int_S f(x)\,\gamma(\cdot)(dx)
  \qquad\text{a.s., for each } f\in C_b(S).
\end{equation*}
Moreover, $\gamma$ can be taken $\sigma(Z_1,Z_2,\dots)$-measurable.
\end{lemma}
We remark that the last sentence of the previous lemma is from the first line of its proof in \cite{berti2004limit}, not from its statement.

\begin{lemma}[L\'evy's upward theorem]\label{lem:levy}
Let $(\mathcal{G}_n)_{n\geq0}$ be a filtration, $\mathcal{G}_\infty=\sigma\big(\bigcup_{n}\mathcal{G}_n\big)$,
and let $Y\in L^1$. Then
\begin{equation*}
  \mathbb{E}[Y\mid\mathcal{G}_n]\to\mathbb{E}[Y\mid\mathcal{G}_\infty]\qquad\text{a.s.\ and in } L^1 .
\end{equation*}
In particular, if $Y$ is $\mathcal{G}_\infty$-measurable, then $\mathbb{E}[Y\mid\mathcal{G}_n]\to Y$ a.s.\ and in $L^1$.
\end{lemma}

The following lemma is from \cite{hunt66}. For a reference in English, see also \cite{blackwell1962merging} Theorem 2.

\begin{lemma}[Hunt's lemma]\label{lem:hunt}
Let $(\mathcal G_n)_{n\ge0}$ be a filtration, $\mathcal{G}_\infty=\sigma\big(\bigcup_{n}\mathcal{G}_n\big)$, and let $(Y_n)_{n\geq 0}$ be random
variables such that $Y_n\to Y$ a.s.\ and $\sup_n|Y_n|\in L^1$. Then
$\mathbb{E}[Y_n\mid\mathcal{G}_n]\to\mathbb{E}[Y\mid\mathcal{G}_\infty]$ a.s.
\end{lemma}

For the following Lemma, see e.g. \cite{DM78}, Chapter~I, Theorem~21.

\begin{lemma}[Functional monotone class theorem]\label{lem:fmct}
Let $S$ be a set and $\mathcal{H}$ a vector space of bounded real functions on $S$ such that
\begin{enumerate}
  \item[(a)] $\mathcal{H}$ contains the constant functions;
  \item[(b)] if $h_m\in\mathcal{H}$, $0\leq h_m\uparrow h$ and $h$ is bounded, then $h\in\mathcal{H}$.
\end{enumerate}
If $\mathcal{C}\subset\mathcal{H}$ is closed under multiplication, then $\mathcal{H}$ contains every
bounded $\sigma(\mathcal{C})$-measurable function.
\end{lemma}

\begin{corollary}\label{cor:fmct}
Let $Q,Q'$ be probability measures on $(S,\mathcal{S})$ and let $\mathcal{C}$ be a class of bounded
$\mathcal{S}$-measurable functions, closed under multiplication, with $\sigma(\mathcal{C})=\mathcal{S}$.
If $\int h\,\di Q=\int h\,\di Q'$ for all $h\in\mathcal{C}$, then $Q=Q'$.
\end{corollary} 

\begin{proof}
The class $\mathcal{H}=\{h \text{ bounded } \mathcal{S}\text{-measurable}:\int h\,\di Q=\int h\,\di Q'\}$ is a
vector space containing $\mathcal{C}$ and the constants (since $Q(S)=Q'(S)=1$), and it is closed under
bounded monotone limits by monotone convergence. By Lemma~\ref{lem:fmct}, $\mathcal{H}$ contains
$\mathbb{I}_A$ for every $A\in\mathcal{S}$.
\end{proof}

\begin{lemma}[Lemma 4.1 of \cite{crimaldi2016fluctuation}] \label{lem:crimaldi_4.1}
Let $\mathcal{G}$ be a filtration and $(X_{n})_{n\geq 0}$ a $\mathcal{G}$-adapted sequences of real r.v.s such that $\mathbb{E}_{n-1}(X_{n})\to X$ a.s. for some real r.v. $X$. Moreover, let $(a_{n})_{n\geq 0}$ and $(b_{n})_{n\geq 0}$ be two sequences of strictly positive real numbers such that $b_{n}\uparrow \infty$ and $\sum_{n=1}^{\infty}\frac{\mathbb{E}[X_{n}^{2}]}{a_{n}^{2}b_{n}^{2}}<+\infty$. Then, we have:
\begin{itemize}
    \item[(a)] If $\frac{1}{b_{n}}\sum_{i=1}^{n} \frac{1}{a_{i}}\to \gamma$ for some constant $\gamma$, then $\frac{1}{b_{n}}\sum_{i=1}^{n}\frac{X_{i}}{a_{i}}\stackrel{a.s.}{\to}\gamma X$.
    \item[(b)] If $\frac{1}{b_{n}}\sum_{i\geq n} \frac{1}{a_{i} b_{i}^{2}}\to \gamma$ for some constant $\gamma$, then $\frac{1}{b_{n}}\sum_{i\geq n}\frac{X_{i}}{a_{i}b_{i}^{2}}\stackrel{a.s.}{\to}\gamma X$.
\end{itemize}
\end{lemma}

\subsection{Almost Supermartingales} \label{sec:almost.supermartingales}
For completeness, we recall some results about almost supermartingales from \cite{robbins1971convergence}. We present these results in a slightly simplified formulation, which is sufficient for the purposes of this paper. Let $(Z_{n})_{n\geq 0}$ and $(\xi_{n})_{n\geq 0}$ be two sequences of random variables adapted to a filtration $(\mathcal{G}_{n})_{n\geq 0}$. Also, let $\xi_{n}$ be non-negative for all $n\geq 0$. Then, $(Z_{n})_{n\geq 0}$ is a \textit{almost $\mathcal{G}$-supermartingale} if for all $n\geq 0$,
\begin{equation}\label{eq:almost.super}
    \mathbb{E}_{n}(Z_{n+1}) \leq Z_{n} + \xi_{n} \ \ \ \ \ \text{a.s.}
\end{equation}

We say $(Z_{n})_{n\geq 0}$ is a \textit{almost $\mathcal{G}$-submartingale} if \eqref{eq:almost.super} holds with reversed inequality, i.e. $\mathbb{E}_{n}(Z_{n+1}) \geq Z_{n} - \xi_{n}$, for all $n\geq 0$. We say that a sequence is an almost $\mathcal{G}$-martingale if it is an almost super and submartingale. 

\cite{robbins1971convergence} have shown that almost supermartingales satisfy a convergence theorem analogous to Doob’s convergence theorem for supermartingales. 

\begin{theorem}[Theorem 1, \cite{robbins1971convergence}] \label{th:robsig}
    If $(Z_{n})_{n\geq 0}$ is a almost $\mathcal{G}$-supermartingale with parameters $(\xi_{n})_{n\geq 0}$, then $\lim_{n\to \infty} Z_{n}$ exists and is finite on the set $\{ \omega\in \Omega: \sum_{n=0}^{\infty} \xi_n(\omega) < \infty \}$.
\end{theorem}

\section{Additional proofs from section~\ref{sec:acid}}

The next lemma shows the equivalence between the two definitions of $\xi$-a.c.i.d. r.v.s stated in Section~\ref{sec:acid}.
\begin{lemma}
Definition~\ref{def:acid.pred} is equivalent to Definition~\ref{def:acid.law}.
\end{lemma}
\begin{proof}
Definition~\ref{def:acid.law} clearly implies Definition~\ref{def:acid.pred}. We show that the converse is also true.
Let $\mathcal{A}_{0}$ the a countable class generating $\mathcal{E}$, which exists because $(E,\mathcal{E})$ is a Polish space, and let $\mathcal{A}$ denotes the algebra generated by $\mathcal{A}_{0}$. Because $\mathcal{A}_{0}$ is countable, $\mathcal{A}$ is also countable. Also, $\mathcal{A}$ generates $\mathcal{E}$. \\
For each $A\in \mathcal{A}$, let us denote by $N_{A}$ the null set over which condition \eqref{eq:acid.pred} fails, i.e. $N_{A}^{c}=\{\omega \in \Omega : \mathbb{P}(X_{n+2}\in A|\mathcal{G}_{n})(\omega) \leq \mathbb{P}(X_{n+1}\in A|\mathcal{G}_{n})(\omega)+ \xi_{n}(\omega)\}$. From Definition~\ref{def:acid.pred}, we know that $\mathbb{P}(N_{A}^{c})=1$ for all $A\in\mathcal{E}$. 
Being $\mathcal{A}$ countable, $N=\cup_{A\in \mathcal{A}}N_{A}$ is also a null set. Moreover, because $\mathcal{A}$ is closed under complements, if $\omega \in N^{c}$, then $\omega \in (N_{A}\cup N_{A^{c}})^{c} = N^{c}_{A}\cap N_{A^{c}}^{c}$ for every $A \in\mathcal{A}$, and therefore $|\mathbb{P}(X_{n+2}\in A|\mathcal{G}_{n})(\omega) - \mathbb{P}(X_{n+1}\in A|\mathcal{G}_{n})(\omega)|\leq \xi_{n}(\omega)
$ for all $A\in \mathcal{A}$. Hence, for all $\omega \in N^{c}$, we have that
\begin{equation*}
    \sup_{A\in \mathcal{A}}|\mathbb{P}(X_{n+2}\in A|\mathcal{G}_{n})(\omega) -  \mathbb{P}(X_{n+1}\in A|\mathcal{G}_{n})(\omega) |\leq  \xi_{n}(\omega)
\end{equation*}
To conclude, we need to show that 
\begin{align} \label{eq:formula100}
    &\sup_{A\in \mathcal{A}}|\mathbb{P}(X_{n+2}\in A|\mathcal{G}_{n})(\omega) -  \mathbb{P}(X_{n+1}\in A|\mathcal{G}_{n})(\omega) | =   \\  &\hspace{2cm} = \sup_{A\in \mathcal{E}}|\mathbb{P}(X_{n+2}\in A|\mathcal{G}_{n})(\omega) -  \mathbb{P}(X_{n+1}\in A|\mathcal{G}_{n})(\omega) |. \nonumber
\end{align}
The ``$\leq$'' inequality is obvious. In order to show the ``$\geq$'' inequality, let us fix $\epsilon>0$. By definition of supremum, there exists $B\in\mathcal{E}$ s.t. 
\begin{align} \label{eq:formula101}
    & \sup_{A\in \mathcal{E}}|\mathbb{P}(X_{n+2}\in A|\mathcal{G}_{n})(\omega) -  \mathbb{P}(X_{n+1}\in A|\mathcal{G}_{n})(\omega) | = \nonumber \\
    & \hspace{2cm} = |\mathbb{P}(X_{n+2}\in B|\mathcal{G}_{n})(\omega) -  \mathbb{P}(X_{n+1}\in B|\mathcal{G}_{n})(\omega) | + \frac{\epsilon}{2}.
\end{align}
Because $\mathcal{A}$ is an algebra generating $\mathcal{E}$, for every $B\in\mathcal{E}$, there exists $B_{0} \in \mathcal{A}$ s.t. $|\mathbb{P}(X_{n+2}\in (B\triangle B_{0})|\mathcal{G}_{n})(\omega) -  \mathbb{P}(X_{n+1}\in (B\triangle B_{0})|\mathcal{G}_{n})(\omega) |\leq \frac{\epsilon}{2}$. Therefore, 
\eqref{eq:formula101} can be upper bounded by
\begin{align*}
    & \leq |\mathbb{P}(X_{n+2}\in B_{0}|\mathcal{G}_{n})(\omega) -  \mathbb{P}(X_{n+1}\in B_{0}|\mathcal{G}_{n})(\omega) | +|\mathbb{P}(X_{n+2}\in (B\triangle B_{0})|\mathcal{G}_{n})(\omega) + \\
    & \hspace{6cm} -  \mathbb{P}(X_{n+1}\in (B\triangle B_{0})|\mathcal{G}_{n})(\omega) |+ \frac{\epsilon}{2}\\
    & \leq |\mathbb{P}(X_{n+2}\in B_{0}|\mathcal{G}_{n})(\omega) -  \mathbb{P}(X_{n+1}\in B_{0}|\mathcal{G}_{n})(\omega) | + \epsilon
\end{align*}
Taking the sup of $B_{0}\in\mathcal{A}$, and because $\epsilon$ is arbitrary, we obtain \eqref{eq:formula100}.
\end{proof}

\subsection{Proof of Example~\ref{ex:garelli}} \label{sec:proof.garelli}
In this subsection, we prove that the sequence of predictive distributions described in Example~\ref{ex:garelli} is $\xi$-a.c.i.d. with $\xi_{n}=\frac{C_{p}}{(n+1)^{2}}$, with $C_{p}$ a known constant, given in the proof of the following proposition.

\begin{proposition}\label{prop:simple}
Let $(X_{n})_{n\geq 1}$ be a sequence of r.v.s with predictive distributions described in Example~\ref{ex:garelli}. Then, there exists a constant $C_p < \infty$, depending only on $p$, such that for every $n \ge 1$,
\[
\text{TV}\big(\mathcal{L}(X_{n+1}\mid X_{1:n}),\ \mathcal{L}(X_{n+2}\mid X_{1:n})\big)
\le \frac{C_{p}}{(n+1)^2}.
\]
\end{proposition}
\begin{proof}
We first prove the bound for $n \ge 2$, and treat the case $n = 1$ at the end.\\
First, we notice that the mean and covariance updates can be written recursively as
\begin{align*}
M_{n+1}& = M_n + \frac{1}{n+1} (X_{n+1} - M_n), \\
Q_{n+1}& = \Big(1-\frac{1}{n+1}\Big)\,Q_n + \frac{1}{n+1}\Big(1-\frac{1}{n+1}\Big) (X_{n+1} - M_n) (X_{n+1} - M_n)^T.
\end{align*}
Write $(X_{n+1} - M_n) = Q_n^{1/2} Z$ with $Z \sim N(0, I_p)$. The total variation distance is invariant under the bijective
affine map $x \mapsto Q_n^{-1/2}(x - M_n)$, hence, conditionally on $X_{1:n}$,
\begin{equation} \label{eq:TV}
\text{TV}\big(\mathcal{L}(X_{n+1}\mid X_{1:n}),\ \mathcal{L}(X_{n+2}\mid X_{1:n})\big)
= \text{TV}\Big(N(0, I_p),\ \E_{Z} \left[ N\left(\frac{1}{n+1} Z, A_Z \right) \right]\Big),
\end{equation}
where
\[A_Z = \left(1-\frac{1}{n+1}\right) I_p + \frac{1}{n+1}\left(1-\frac{1}{n+1}\right) Z Z^T.\]
In particular, the right-hand side depends only on $(n, p)$ and not on the data $X_{1:n}$. We need to upper bound the TV distance between a standard Normal and a mixture of Normal distributions (with standard Normal mixing measure).
\\
To lighten the formulas, we write $\lambda = \frac{1}{n+1} \in (0, \frac13]$, and 
denote by $\varphi$ the density of the standard Normal $N(0, I_p)$, and by $\varphi_{Z,\lambda}$ the density of the Normal distribution $N(\lambda Z, A_Z)$. Note that the TV in equation \eqref{eq:TV} can be written as
\begin{equation}\label{eq:TV0}
    \text{TV} =\frac{1}{2}\E_{\varphi}\Big[ \Big| \frac{\E_{Z}[\varphi_{Z,\lambda}]}{\varphi} -1 \Big| \Big].
\end{equation}
We start by noticing that if $Y$ has the mixture law in \eqref{eq:TV}, i.e.\ $Y \mid Z \sim N(\lambda Z, A_Z)$, then, for every Borel set $B$,
by Tonelli's theorem,
\[
\mathbb{P}(Y \in B) = \E_Z\Big[ \int_B \varphi_{Z,\lambda}(y)\, dy \Big] = \int_B \E_Z\big[ \varphi_{Z,\lambda}(y) \big] dy ,
\]
so $Y$ has density $f(y) = \E_Z[\varphi_{Z,\lambda}(y)]$. Since $\varphi > 0$ and $\varphi_{Z,\lambda} > 0$, denoting by $L_Z := \log(\varphi_{Z,\lambda}/\varphi)$ the log-likelihood ratio, we have
$\varphi_{Z,\lambda}(y) = \varphi(y)\, e^{L_Z(y)}$, and, as $\varphi(y)$ does not depend on $Z$,
\[
f(y) = \E_Z[\varphi_{Z,\lambda}(y)] = \varphi(y)\, \E_Z\big[ e^{L_Z(y)} \big].
\]
Using this fact in \eqref{eq:TV0},
\begin{align}\label{eq:bound00}
\text{TV} & 
=\frac{1}{2} \E_\varphi \Big[ \big| \E_Z\big[e^{L_Z} - 1\big] \big| \Big].
\end{align}
We compute the log-likelihood ratio $L_Z$. Plugging in the multivariate Normal densities $\varphi$ and $\varphi_{Z,\lambda}$, and using the fact that, for $A_Z = (1-\lambda)\big(I_p + \lambda ZZ^T\big)$, 
\[
\det A_Z = (1-\lambda)^p (1 + \lambda |Z|^2), \quad
A_Z^{-1} = \frac{1}{1-\lambda}\Big( I_p - \frac{\lambda}{1+\lambda |Z|^2} ZZ^T \Big), \quad
A_Z^{-1} Z = \frac{Z}{(1-\lambda)(1+\lambda |Z|^2)},
\]
the log-likelihood ratio $L_Z$ is the quadratic polynomial in $y$
\begin{align*}
L_Z(y) &= -\tfrac12 \log\det A_Z - \tfrac12 (y - \lambda Z)^T A_Z^{-1} (y - \lambda Z) + \tfrac12 |y|^2 \\
&= -\frac{\lambda\, |y|^2}{2(1-\lambda)}
+ \frac{\lambda\, (Z\cdot y)^2}{2(1-\lambda)(1+\lambda |Z|^2)}
+ \frac{\lambda\, Z\cdot y}{(1-\lambda)(1+\lambda |Z|^2)}
- \frac{\lambda^2 |Z|^2}{2(1-\lambda)(1+\lambda |Z|^2)}\\
& \hspace{1cm}- \frac{p}{2}\log(1-\lambda) - \frac12 \log(1+\lambda |Z|^2).
\end{align*}
We denote by $D_Z$ the part of $L_Z$ that is linear in $\lambda$:
\[
D_Z(y) = \lambda \Big[ \frac{p - |Z|^2}{2} + Z\cdot y + \frac12\big( (Z\cdot y)^2 - |y|^2 \big) \Big].
\]
For every fixed $y$, since $\E_Z[|Z|^2] = p$, $\E_Z[Z] = 0$ and $\E_Z\big[(Z\cdot y)^2\big] = |y|^2$, we have $\E_Z [ D_Z(y) ]= 0$. Using this last fact in \eqref{eq:bound00}, and using Jensen's inequality and Tonelli's theorem for the inequality, we obtain
\begin{align}\label{eq:bound0}
\text{TV} &  
=\tfrac12\, \E_\varphi \Big[ \big| \E_Z\big[e^{L_Z} - 1\big] \big| \Big]
= \tfrac12\, \E_\varphi \Big[ \big| \E_Z\big[e^{L_Z} - 1 - D_Z\big] \big| \Big] \nonumber \\
& \le \tfrac12\, \E_Z \Big[ \E_\varphi \big[\big| e^{L_Z} - 1 - D_Z \big|\big] \Big].
\end{align}
We aim to upper bound the inner expectation in \eqref{eq:bound0}. For real $a, b$ we have $|e^a - 1 - a| \le \frac12 a^2 \max(1, e^a)$, hence
$|e^a - 1 - b| \le |a - b| + \tfrac12 a^2 (1 + e^a)$. Using this inequality with $a = L_Z(y)$, $b = D_Z(y)$, and using $\E_\varphi\big[L_Z^2\, e^{L_Z}\big] = \E_{\varphi_{Z,\lambda}}\big[L_Z^2\big]$,
\begin{equation}\label{eq:bound}
\E_\varphi \left[ \big| e^{L_Z} - 1 - D_Z \big| \right]
\le \E_\varphi \left[  |L_Z - D_Z| \right] + \tfrac12\, \E_\varphi \left[ L_Z^2 \right] + \tfrac12 \E_{\varphi_{Z,\lambda}} \left[ L_Z^2 \right].
\end{equation}
We upper bound the three expectations on the right-hand side of \eqref{eq:bound}.
First, let us consider $\E_\varphi \left[  |L_Z - D_Z| \right]$. Subtracting $D_Z$ from $L_Z$ coefficient by coefficient,
\begin{align*}
L_Z - D_Z &= -\frac{\lambda^2}{2(1-\lambda)}\, |y|^2
+ \frac{\lambda^2}{2(1-\lambda)}\Big(1 - \frac{|Z|^2}{1+\lambda |Z|^2}\Big) (Z\cdot y)^2
+ \frac{\lambda^2}{1-\lambda}\Big(1 - \frac{|Z|^2}{1+\lambda |Z|^2}\Big)\, Z\cdot y \\
&\quad - \frac{\lambda^2 |Z|^2}{2(1-\lambda)(1+\lambda |Z|^2)} + \frac p2\big(-\log(1-\lambda) - \lambda\big)
+ \frac{1}{2}\big(\lambda |Z|^2 - \log(1+\lambda |Z|^2)\big).
\end{align*}
Using $\lambda \le \frac13$ (so $\frac{1}{1-\lambda} \le \frac32$), $\frac{1}{1+\lambda |Z|^2} \le 1$,
$\big|1 - \frac{|Z|^2}{1+\lambda |Z|^2}\big| \le 1 + |Z|^2$, $0 \le -\log(1-\lambda) - \lambda \le \frac{\lambda^2}{2(1-\lambda)}$ and
$0 \le \lambda |Z|^2 - \log(1+\lambda |Z|^2) \le \frac{\lambda^2 |Z|^4}{2}$, every coefficient is $\lambda^2$ times a quantity bounded by a
constant multiple of $(1+|Z|^2)^2$. Since moreover $(Z\cdot y)^2 \le |Z|^2\,|y|^2$ and $|Z\cdot y| \le \frac12 (|Z|^2 + |y|^2) \le (1+|Z|^2)(1+|y|^2)$,
we obtain
\begin{equation}\label{eq:bound_difference}
|L_Z - D_Z| \le (p+4) \lambda^2 (1+|Z|^2)^2 \big(1 + |y|^2\big).
\end{equation}
Under a standard Normal $\varphi$, $|y|^2 \sim \chi^2_p$, and $\E_\varphi\big[1 + |y|^2\big] = p + 1$. Hence,
\begin{equation}\label{eq:bound_difference1}
\E_\varphi \left[  |L_Z - D_Z| \right] \le (p+4) \lambda^2 (1+|Z|^2)^2  \E_\varphi \left[  1 + |y|^2 \right] = (p+4) (p+1) \lambda^2 (1+|Z|^2)^2 .
\end{equation}
Regarding the other two expectations in \eqref{eq:bound}, we first notice that each coefficient of $L_Z$ is $\lambda$ times a quantity bounded by a constant multiple of $1+|Z|^2$
(using also $-\frac p2 \log(1-\lambda) \le \frac{3\lambda p}{4}$ and $\frac12 \log(1+\lambda |Z|^2) \le \frac{\lambda |Z|^2}{2}$), so
\begin{equation}\label{eq:bound_L}
|L_Z| \le (p+3) \lambda (1+|Z|^2) \big(1 + |y|^2\big).
\end{equation}
Again, under $\varphi$, $|y|^2 \sim \chi^2_p$, and, from the moments of $\chi^2_p$,
$\E_\varphi\big[(1 + |y|^2)^2\big] = 1 + 2p + p(p+2) \le (p+2)^2$. Hence, 
\begin{equation}\label{eq:bound_L_phi}
\mathbb{E}_{\varphi}[L_Z^{2}] \le (p+3)^{2} \lambda^2 (1+|Z|^2)^{2} \mathbb{E}_{\varphi}\big[(1 + |y|^2)^{2}\big]\leq (p+3)^{2} (p+2)^{2} \lambda^2 (1+|Z|^2)^{2}.
\end{equation}
Under the density $\varphi_{Z,\lambda}$, we can write $Y = \lambda Z + A_Z^{1/2} W$ with $W \sim N(0, I_p)$. The largest eigenvalue of $A_Z$ is
$(1-\lambda)(1+\lambda |Z|^2) \le 1 + |Z|^2$, so $|Y|^2 \le 2\lambda^2 |Z|^2 + 2(1+|Z|^2)|W|^2 \le 2(1+|Z|^2)\big(1 + |W|^2\big)$, hence
$1 + |Y|^2 \le 3(1+|Z|^2)\big(1 + |W|^2\big)$ and
\[
\E_{\varphi_{Z,\lambda}}\big[(1 + |Y|^2)^2\big] \le 9 (1+|Z|^2)^2\, \E\big[(1 + |W|^2)^2\big] \le 9 (p+2)^2 (1+|Z|^2)^2 .
\]
Therefore, using this bound and \eqref{eq:bound_L},
\begin{equation}\label{eq:bound_L_phiZ}
\E_{\varphi_{Z,\lambda}}[ L_Z^2] \le (p+3)^{2} \lambda^2 (1+|Z|^2)^{2} \E_{\varphi_{Z,\lambda}}\big[ (1 + |Y|^2)^{2}\big]\leq 9 (p+3)^{2} (p+2)^2 \lambda^2 (1+|Z|^2)^{4}.
\end{equation}
Therefore, using the bounds \eqref{eq:bound_difference1}, \eqref{eq:bound_L_phi} and \eqref{eq:bound_L_phiZ} in \eqref{eq:bound}, we obtain
\begin{align} \label{eq:bound4}
\E_\varphi\big[|e^{L_Z} - 1 - D_Z|\big]
& \le \lambda^2 \Big[ (p+4)(p+1)(1+|Z|^2)^2 + \frac{1}{2} (p+3)^2 (p+2)^2 (1+|Z|^2)^2 \nonumber\\
& \hspace{2cm} +\frac{9}{2} (p+3)^2 (p+2)^2 (1+|Z|^2)^4 \Big] \nonumber\\
&\le c_p\, \lambda^2 (1+|Z|^2)^4,
\end{align}
with $c_p = (p+4)(p+1) + 5(p+3)^2(p+2)^2$. Finally, using \eqref{eq:bound4} in \eqref{eq:bound0}, we obtain 
\[
\text{TV} \leq \frac{1}{2}\, \E_Z\Big[ \E_\varphi\big[|e^{L_Z} - 1 - D_Z|\big] \Big]
\le \frac{c_p\, \E_Z\big[(1+|Z|^2)^4\big]}{2} \lambda^2 = \frac{C_p}{(n+1)^2},
\]
with $C_p := \tfrac12\, c_p\, \E_Z\big[(1+|Z|^2)^4\big] < \infty$, since $|Z|^2 \sim \chi^2_p$ has finite moments of all orders. This proves the bound for $n \geq 2$.
\\
Finally, consider $n = 1$. Since $c_p$ is increasing in $p$, $c_p \geq c_1 = 730$, and $\E_Z\big[(1+|Z|^2)^4\big] \ge 1$, so
$C_p \ge 365$. As the total variation distance is always at most $1$, $\text{TV}\big(\mathcal{L}(X_{2}\mid X_{1}), \mathcal{L}(X_{3}\mid X_{1})\big) \le 1 \leq C_p/4$.
\end{proof}

\section{Additional proofs from Section~\ref{sec:slln.ae} and Section~\ref{sec:abs.cont}}

The following two lemmas are used in the proof of Theorem~\ref{thm:asym.exch}. 
\begin{lemma}\label{lem:gen.class}
Let $E$ be a Polish space and let $D_0\subset C_b(E)$ be a countable convergence-determining class
containing the constant $1$ and closed under multiplication. Then, for every $k\geq1$, the class
\begin{equation*}
  \mathcal{C}_k=\Big\{x\mapsto\prod_{j=1}^k f_j(x_j):\ f_1,\dots,f_k\in D_0\Big\}
\end{equation*}
is closed under multiplication, contains $1$, and satisfies $\sigma(\mathcal{C}_k)=\mathcal{B}(E^k)$.
\end{lemma}
\begin{proof}
The first two properties are inherited from $D_0$. We first show that $D_0$ separates the points
of $E$: if $f(x)=f(y)$ for all $f\in D_0$, then the constant sequence $\delta_x$ converges to $\delta_y$
along $D_0$, hence weakly, because $D_0$ is convergence-determining, so $\delta_x=\delta_y$ and $x=y$.

Next, we show that $\sigma(D_0)=\mathcal{B}(E)$. The inclusion $\subseteq$ holds since each $f\in D_0$
is continuous. Conversely, write $D_0=\{f_1,f_2,\dots\}$ and let $\Phi:E\to\mathbb{R}^{\mathbb{N}}$,
$\Phi(x)=(f_m(x))_{m\geq1}$. Then $\Phi$ is continuous, injective since $D_0$ separates points, and
$\sigma(D_0)=\Phi^{-1}(\mathcal{B}(\mathbb{R}^{\mathbb{N}}))$. For $A\in\mathcal{B}(E)$, the set $\Phi(A)$ is
Borel by the Lusin--Souslin theorem (\cite{Kec95}, Theorem~15.1), so, by injectivity,
$A=\Phi^{-1}(\Phi(A))\in\sigma(D_0)$.

Finally, each element of $\mathcal{C}_k$ is continuous, so $\sigma(\mathcal{C}_k)\subseteq\mathcal{B}(E^k)$.
Conversely, taking $f_i\equiv1$ for $i\neq j$, the class $\mathcal{C}_k$ contains $x\mapsto f(x_j)$ for every
$f\in D_0$ and $j\leq k$. Hence $\sigma(\mathcal{C}_k)\supseteq\sigma(D_0)^{\otimes k}=\mathcal{B}(E)^{\otimes k}
=\mathcal{B}(E^k)$, where the last equality holds because $E$ is separable, and $^{\otimes k}$ denotes the $k$-fold product $\sigma$-algebra.
\end{proof}

\begin{lemma}\label{lem:ae.prod}
Let $(X_n)_{n\geq1}$ be a $\xi$-a.c.i.d.\ sequence with
$\sum_{n\geq0}\xi_n<\infty$ a.s. For $f\in C_b(E)$, let $V_f:=\lim_n\int_E f\,dP_n$, which exists a.s.\
by the first part of the proof of Theorem~\ref{thm:slln} and is $\mathcal{G}_\infty$-measurable as an
a.s. limit of $\mathcal{G}_n$-measurable random variables. Then, for all $k\geq1$ and
$f_1,\dots,f_k\in C_b(E)$,
\begin{equation}\label{eq:ae.prod}
  \mathbb{E}\Big[\prod_{j=1}^k f_j(X_{n+j})\,\Big|\,\mathcal{G}_n\Big]\to\prod_{j=1}^k V_{f_j}
  \qquad\text{a.s.}
\end{equation}
\end{lemma}
\begin{proof}
We argue by induction on $k$. The case $k=1$, since $\mathbb E[f_1(X_{n+1})\mid\mathcal G_n]=\int_E f_1\,dP_n$, is immediate. Let $k \geq 2$ and assume
\eqref{eq:ae.prod} for $k-1$. Conditioning on $\mathcal{G}_{n+k-1}$, since $f_{j}(X_{n+j})$ is $\mathcal{G}_{n+k-1}$-measurable for $j<k$,
\begin{align*} \mathbb{E}\Big[\prod_{j=1}^k f_j(X_{n+j})\Big|\mathcal{G}_n\Big] &=   \mathbb{E}\Big[ \mathbb{E} \big[ \prod_{j=1}^k f_j(X_{n+j})\big| \mathcal{G}_{n+k-1} \big]   \Big|\mathcal{G}_n\Big]
  \\
&=\mathbb{E}\Big[\prod_{j=1}^{k-1} f_j(X_{n+j})\int_E f_k\,dP_{n+k-1}\Big|\mathcal{G}_n\Big].
\end{align*}
Write $\int_E f_k dP_{n+k-1}=\mathbb{E}[V_{f_k}\mid\mathcal{G}_n]+R_n$, where $R_{n}=\int_{E}f_{k}dP_{n+k-1}- \mathbb{E}[V_{f_{k}}|\mathcal{G}_{n}]$. By the definition of $V_{f_k}$ and Lemma~\ref{lem:levy} ($V_{f_k}$ is $\mathcal{G}_\infty$-measurable), $R_n\to V_{f_k}-V_{f_k}=0$ a.s., and
$|R_n|\leq2\lVert f_k\rVert_\infty$. Hence,
\begin{align*}
  \mathbb{E}\Big[\prod_{j=1}^k f_j(X_{n+j})\,\Big|\,\mathcal{G}_n\Big]
  =\mathbb{E}[V_{f_k}\mid\mathcal{G}_n]\,
   \mathbb{E}\Big[\prod_{j=1}^{k-1} f_j(X_{n+j})\,\Big|\,\mathcal{G}_n\Big]+\Delta_n, 
\end{align*}
with $|\Delta_n|\leq\prod_{j=1}^{k-1}\lVert f_j\rVert_\infty\,\mathbb{E}\big[|R_n|\,\big|\,\mathcal{G}_n\big]$. The first term converges a.s.\ to $V_{f_k}\prod_{j<k}V_{f_j}$ by Lemma~\ref{lem:levy} and the
inductive assumption, while $\Delta_n\to0$ a.s.\ by Lemma~\ref{lem:hunt} applied to $Y_n=|R_n|$.
This proves \eqref{eq:ae.prod}.
\end{proof}

The following lemma is used in Example~\ref{ex:exch_rvs}.

\begin{lemma}\label{lem:newlemma}
    Let $(X_{n})_{n\geq 1}$ be a sequence of r.v.s such that $\frac{1}{n}\sum_{i=1}^{n}X_{i} \stackrel{\text{a.s.}}{\to} \mu$, $(\sigma^{2}_{i})_{i\in \mathbb{N}}$ a sequence of constants converging to a value $\sigma^{2}_{*}$, and $(\sigma_{\mu}^{2},\sigma_{r}^{2})$ positive constants. Then,
    \begin{equation*}
        \left(\frac{1}{\sigma_{\mu}^{2} }
        + \sum_{i=1}^{n} \frac{1}{\sigma_{i}^{2}+\sigma_{r}^{2} }\right)^{-1}\sum_{i=1}^{n} \frac{X_{i}}{\sigma^{2}_{i}+\sigma^{2}_{r}} \stackrel{\text{a.s.}}{\to} \mu
    \end{equation*}
\end{lemma}
\begin{proof}
    Denote $w_{i}=\frac{1}{\sigma^{2}_{i}+\sigma^{2}_{r}}$, $W_{n}=\sum_{i=1}^{n}w_{i}$ and $c=1/\sigma_{\mu}^{2}$. \\
    We need to prove that $\frac{1}{c+W_{n}}\sum_{i=1}^{n}w_{i} X_{i} \stackrel{\text{a.s.}}{\to} \mu$. From 
    \begin{equation*}
   \frac{1}{c+W_{n}}\sum_{i=1}^{n}w_{i} X_{i} =\frac{W_{n}}{c+W_{n}}\left( \mu +\frac{1}{W_{n}} \sum_{i=1}^{n}w_{i} (X_{i}-\mu )\right),
    \end{equation*}
    because $W_{n}\to \infty$ (from $\sigma_{r}^{2}>0$), it is enough to show that $\frac{1}{W_{n}} \sum_{i=1}^{n}w_{i} (X_{i}-\mu ) \stackrel{\text{a.s.}}{\to} 0$. \\
    Let $S_{n} =\sum_{i=1}^{n}(X_{i}-\mu)$. Summing by parts,
    \begin{equation} \label{eq:formula.newlemma}
     \frac{1}{W_{n}} \sum_{i=1}^{n}w_{i} (X_{i}-\mu ) = \frac{1}{W_{n}} w_{n} S_{n} + \frac{1}{W_{n}} \sum_{i=1}^{n-1}S_{i}(w_{i}-w_{i+1})
    \end{equation}
    The first term goes a.s. to zero. Indeed, 
    because $\sigma_{i}^{2} \to \sigma_{*}^{2}$, $w_{i}\to w^{*}=\frac{1}{\sigma^{2}_{*}+\sigma^{2}_{r}}$, and $0<w^{*}<\infty$. So, there are constants $m,M$ such that $m\leq w_{n}\leq M$  and $m n\leq W_{n} \leq M n$, for $n$ large enough. Together these inequalities imply, $\frac{w_{n}}{W_{n}}\leq \frac{M}{mn}$, for $n$ large. Therefore, $ \frac{w_{n}}{W_{n}}  S_{n} \leq \frac{M}{m} \frac{S_{n}}{n}$, and $\frac{S_{n}}{n} \stackrel{\text{a.s.}}{\to} 0$ by assumption.\\
    The second term goes a.s. to zero. Fix $\epsilon >0$ arbitrary. Because  $\frac{S_{n}}{n} \stackrel{\text{a.s.}}{\to} 0$, $\exists K$ (depending on $\omega$) such that $|S_{n}|\leq n \epsilon$, for all $n\geq K$. Also, because $(w_{i})_{i\in\mathbb{N}}$ converges to $w^{*}$, the partial sums of $|w_{i}-w_{i+1}|$ are bounded (by a constant $C$). Then, for $n$ large (using $\frac{n}{W_{n}}\leq \frac{1}{m}$)
    \begin{equation*}
        \frac{1}{W_{n}} |\sum_{i=1}^{n-1}S_{i}(w_{i}-w_{i+1}) |\leq  \frac{1}{W_{n}} C\epsilon n \leq \frac{C}{m} \epsilon
    \end{equation*}
    Because $\epsilon$ is arbitrary the left hand side converges to zero and so the second term in \eqref{eq:formula.newlemma}.
\end{proof}

\section{Central Limit Theorems for Almost Martingales}

In this section, we state two CLTs applicable to almost martingale arrays. These two theorems are used in the proofs of the CLTs for $\xi$-a.c.i.d. sequences, Theorems~\ref{thm:clt_pred_limit}, \ref{thm:clt_empirical_predictive} and \ref{thm:clt_empirical_limit}. Theorem~\ref{thm:HH3.2} corresponds to Theorems 3.2 and 3.3 of \cite{HH80}, in which the martingale assumption, i.e. $\mathbb{E}_{i-1}\left(Y_{n,i}\right)=0$ for all $n\in \mathbb{N}$ and $i\leq k_{n}$, is replaced by the asymptotic conditions \eqref{eq:stableCLT_conv_exp}
and \eqref{eq:stableCLT_conv_exp2}. This kind of generalization mimics Corollary 2.6 of \cite{ML74}, which is stated in terms of convergence in distribution while the theorems considered here are in terms of stable convergence. In the same way, Theorem~\ref{thm:HH3.6} is Theorem 3.6 of \cite{HH80} with the martingale assumption replaced by \eqref{eq:HH3.6_conv_exp}
and \eqref{eq:HH3.6_conv_exp2}, and considering only a single (not doubly) infinite sequence, which is enough for our purposes. We could not find a reference for these versions of CLTs for almost martingales. We include them here with a proof for completeness.

\begin{theorem}\label{thm:HH3.2}
Let $\{S_{n,i},\mathcal{F}_{n,i};n\geq 1,1\leq i\leq k_{n}\}$ be a triangular 
array with differences $Y_{n,i}:=S_{n,i}-S_{n,i-1}$ and let $\eta^{2}$ be an a.s. finite non-negative r.v. Suppose that
\begin{multiequation}
\label{eq:stableCT_nestedsigmafields}
    \mathcal{F}_{n,i} \subseteq \mathcal{F}_{n+1,i},  1\leq i \leq k_n,n \geq 1,
 \hspace{0.8cm}\\
 \label{eq:stableCLT_conv_sum2}
\sum_{i=1}^{k_n}Y_{n,i}^{2} \stackrel{p} {\to}  \eta^{2},
\end{multiequation}
\begin{multiequation}
    \label{eq:stableCT_max_conv}
\max_{i\leq k_{n}}|Y_{n,i}|\stackrel{p}{\to}  0, \hspace{0.8cm}\\ \label{eq:stableCLT_unif_bound}
\sup_{n\geq 1}\mathbb{E}\left(\max_{i\leq k_{n}}Y_{n,i}^{2} \right)<\infty 
\end{multiequation}
\begin{multiequation}
    \label{eq:stableCLT_conv_exp}
\sum_{i=1}^{k_{n}}\mathbb{E}_{i-1}\left(Y_{n,i}\right) \stackrel{p} {\to}  0, \hspace{0.8cm}\\ \label{eq:stableCLT_conv_exp2}
\sum_{i=1}^{k_{n}}\mathbb{E}_{i-1}\left(Y_{n,i}\right)^{2}\stackrel{p} {\to}  0.
\end{multiequation}
Then, $S_{n,k_{n}}=\sum_{i=1}^{k_n} Y_{n,i}\stackrel{d}{\to}Z$ (stably), where the random variable $Z$ has characteristic function $\mathbb{E}\left(\exp(-\frac{1}{2}\eta^{2}t^{2}) \right)$. \\
Moreover, if $\mathbb{P}(\eta^{2}>0)=1$, then
\begin{equation} \label{eq:stableCLT_conclusion}
\left(\sum_{i=1}^{k_{n}}Y_{n,i}^{2}\right)^{-1/2}S_{n,k_{n}} \stackrel{d}{\to}\N(0,1) \ \ \ \text{(mixing)} 
\end{equation}
\end{theorem}
\begin{proof}
Define $\tilde{Y}_{n,i}:=Y_{n,i}- \mathbb{E}_{i-1}\left(Y_{n,i}\right)$. $\tilde{Y}_{n,i}$ is a martingale difference array, and we will check that it satisfies the conditions of Theorem 3.2 of \cite{HH80}. 
\\
Regarding condition (3.18), we need to show that $\max_{i\leq k_{n}}|\tilde{Y}_{n,i}|\stackrel{p}{\to}0$.
    \begin{align*}
       \max_{i\leq k_{n}}|\tilde{Y}_{n,i}| & = \max_{i \leq k_{n}}|Y_{n,i}-\mathbb{E}_{i-1}(Y_{n,i})|  \leq \max_{i\leq k_{n}}|Y_{n,i}| + \max_{i\leq k_{n}}|\mathbb{E}_{i-1}(Y_{n,i})|  \\
        & \leq \max_{i\leq k_{n}}|Y_{n,i}| + \max_{i\leq k_{n}}\mathbb{E}_{i-1}(|Y_{n,i}|)  \leq \max_{i\leq k_{n}}|Y_{n,i}| + \max_{i \leq k_{n}}\mathbb{E}_{i-1}(\max_{j\leq k_{n}}|Y_{n,j}|).
    \end{align*}

The first term converges to zero in probability from assumption~\eqref{eq:stableCT_max_conv}. We need to show that the second term also converges to zero in probability. By noticing that 
$\mathbb{E}_{i}\left( \mathbb{E}_{i-1}\left( \max_{1\leq j\leq k_{n}}|Y_{n,j}| \right) \right) = \mathbb{E}_{i-1}\left( \max_{1\leq j\leq k_{n}}|Y_{n,j}| \right)$, $\left(\mathbb{E}_{i-1}\left( \max_{j\leq k_{n}}|Y_{n,j}|\right)\right)_{i\leq k_{n}}$ is a martingale. Applying Doob's inequality, we obtain,
\begin{align*}
    \mathbb{P}\left( \max_{i\leq k_{n}}\mathbb{E}_{i-1}( \max_{j \leq k_{n}}|Y_{n,j}|) > \epsilon \right) & \leq\frac{1}{\epsilon}\mathbb{E}\left(\mathbb{E}_{k_{n}-1}(\max_{j \leq k_{n}}|Y_{n,j}|)\right)  \leq \frac{1}{\epsilon} \mathbb{E}\left( \max_{j \leq k_{n}}|Y_{n,j}|\right).
\end{align*}
Finally,
\begin{align*}
    \mathbb{E}\left( \max_{j\leq k_{n}}|Y_{n,j}|\right) & \leq \epsilon' + \mathbb{E}\left( \max_{j \leq k_{n}}|Y_{n,j}| \mathbb{I}(\max_{j \leq k_{n}}|Y_{n,j}| > \epsilon' )\right) \\
    & \leq \epsilon' + 
 \left( \mathbb{E}\left( (\max_{j \leq k_{n}}|Y_{n,j}|)^{2}\right)  \mathbb{P}\left(\max_{j\leq k_{n}}|Y_{n,j}| > \epsilon' \right) \right)^{1/2} \\
    & = \epsilon' + 
 \left( \mathbb{E}\left( \max_{j \leq k_{n}}Y_{n,j}^{2}\right)  \mathbb{P}\left(\max_{j \leq k_{n}}|Y_{n,j}| > \epsilon' \right) \right)^{1/2},
\end{align*}
where the last inequality follows from the Cauchy-Schwartz inequality. Finally, the second term converges to zero, because $\mathbb{P}\left(\max_{j\leq k_{n}}|Y_{n,j}| > \epsilon' \right) \stackrel{p}{\to}0$, from assumption~\eqref{eq:stableCT_max_conv}, and\\ $\sup_{n\geq 1} \mathbb{E}\left( \max_{j\leq k_{n}}Y_{n,j}^{2}\right)<\infty$, from  assumption~\eqref{eq:stableCLT_unif_bound}. \\
Regarding condition (3.19) of Theorem 3.2 \cite{HH80}, we need to show that $\sum_{i\leq k_{n}}\tilde{Y}^{2}_{n,i}\stackrel{p}{\to}\eta^{2}$
\begin{align*}
    \sum_{i\leq k_{n}}\tilde{Y}^{2}_{n,i} &  = \sum_{i \leq k_{n}}Y_{n,i}^{2} -2 \sum_{i \leq k_{n}}Y_{n,i}\mathbb{E}_{i}\left( Y_{n,i}\right) + \sum_{i \leq k_{n}}\mathbb{E}_{i}\left(Y_{n,i}\right)^{2}.
\end{align*}
The first term converges in probability to $\eta^{2}$ from assumption \eqref{eq:stableCLT_conv_sum2}, while the last term converges to 0 from assumption \eqref{eq:stableCLT_conv_exp2}. Regarding the middle term, we can upper bound it by
\begin{align*}
    \left| \sum_{i\leq k_{n}}Y_{n,i}\mathbb{E}_{i}\left( Y_{n,i}\right)  \right| & \leq    \sum_{i\leq k_{n}}|Y_{n,i}| \sum_{i \leq k_{n}}|\mathbb{E}_{i}\left( Y_{n,i}\right)  |  \leq  \left(  \sum_{i\leq k_{n}}Y_{n,i} ^{2}\right)^{1/2} \left( \sum_{i \leq k_{n}}\mathbb{E}_{i}\left( Y_{n,i}\right)^{2}\right)^{1/2}, 
\end{align*}
where the last inequality follows from Cauchy-Schwartz. The last term converges in probability to zero, again from assumptions \eqref{eq:stableCLT_conv_sum2} and \eqref{eq:stableCLT_conv_exp2}.  \\
Finally, condition (3.20) of Theorem 3.2 \cite{HH80} is verified in the same way as in Corollary (2.6) of \cite{ML74}.\\
Thus, $\sum_{i\leq k_{n}} \tilde{Y}_{n,i}\stackrel{d}{\to}Z$ (stably), where $Z$ has characteristic function $\mathbb{E}\left(\exp(-\frac{1}{2}\eta^{2}t^{2}) \right)$ by Theorem 3.2 of \cite{HH80}. Slutsky's theorem and Assumption \eqref{eq:stableCLT_conv_exp} imply 
$\sum_{i\leq k_{n}} Y_{n,i}\stackrel{d}{\to}Z$ (stably). Finally, the last part of the Theorem, equation~\eqref{eq:stableCLT_conclusion}, is proven exactly as Theorem 3.3 of \cite{HH80}.
\end{proof}

\begin{theorem}\label{thm:HH3.6}
Let $\{S_{n,i},\mathcal{F}_{n,i};n\geq 1,1\leq i <\infty\}$ be a doubly infinite  
array with differences $Y_{n,i}:=S_{n,i}-S_{n,i-1}$ s.t.
\begin{equation} \label{eq:stableCLT1_sup}
    \sup_{n,i \geq 1}\mathbb{E}(S_{n,i}^{2})<\infty, 
\end{equation}
and let $\eta^{2}$ be an a.s. finite non-negative r.v.. Suppose that, 
\begin{multiequation}
    \label{eq:HH3.6_nestedsigmafields}
    \mathcal{F}_{n,i} \subseteq \mathcal{F}_{n+1,i} \hspace{0.2cm} \text{for} \hspace{0.1cm}  i \geq 1, \ n \geq 1, \hspace{0.8cm}\\
    \label{eq:HH3.6_conv_sum2}
\sum_{i \geq 1}Y_{n,i}^{2} \stackrel{p} {\to}  \eta^{2},
\end{multiequation}
\begin{multiequation}
    \label{eq:HH3.6_max_conv}
\sup_{i\geq 1}|Y_{n,i}|\stackrel{p}{\to}  0, \hspace{0.8cm}\\
\label{eq:HH3.6_unif_bound}
\sup_{n\geq 1}\mathbb{E}\left(\sup_{i\geq 1}Y_{n,i}^{2} \right)<\infty 
\end{multiequation}
\begin{multiequation}
    \label{eq:HH3.6_conv_exp}
\sum_{i=1}^{\infty}\mathbb{E}_{i-1}\left(Y_{n,i}\right) \stackrel{p} {\to}  0, \hspace{0.8cm}\\
\label{eq:HH3.6_conv_exp2}
\sum_{i=1}^{\infty}\mathbb{E}_{i-1}\left(Y_{n,i}\right)^{2}\stackrel{p} {\to}  0.
\end{multiequation}
Then, $S_{n,\infty}=\sum_{i=1}^{\infty}Y_{n,i}\stackrel{d}{\to}Z$, where $Z$ has characteristic function $\mathbb{E}\left(\exp(-\frac{1}{2}\eta^{2}t^{2}) \right)$.\\
Moreover, if $\mathbb{P}(\eta^{2}>0)=1$, then
\begin{equation} \label{eq:HH3.6_conclusion}
\left(\sum_{i=1}^{\infty}Y_{n,i}^{2}\right)^{-1/2}S_{n,\infty} \stackrel{d}{\to}\text{N}(0,1). 
\end{equation}
\end{theorem}
\begin{proof}
    The proof follows that of Theorem 3.6 of \cite{HH80}. In view of \eqref{eq:stableCLT1_sup}, there exists a sequence $(k_{n})_{n\geq1}$ increasing to $\infty$, such that, as $n\to \infty$,
    \begin{equation} \label{eq:HH3.6.proof}
        \mathbb{E}\left[ \sum_{i>k_{n}} Y_{n,i}^{2}\right] \to 0.  
    \end{equation}
We can define the \emph{finite} array $\{S'_{n,i},\mathcal{F}'_{n,i};n\geq 1,1\leq i \leq k_{n}\}$, with $S'_{n,i}=S_{n,i}$ and $\mathcal{F}'_{n,i}=\mathcal{F}_{n,i}$. This array satisfies the assumptions of Theorem~\ref{thm:HH3.2}, and so $S'_{n,k_{n}} \to Z$ (stably), where $Z\sim N(0,\eta^{2})$. Moreover, \eqref{eq:HH3.6.proof} implies
\begin{equation*}
    S_{n,\infty}-S'_{n,k_{n}}=\sum_{i>k_{n}}Y_{n,i} \stackrel{p}{\to} 0.
\end{equation*}
From Slutsky theorem, we have that $S_{n,\infty}\stackrel{d}{\to} Z$. \\
Finally, if $\mathbb{P}(\eta^{2}>0)=1$, from the last part of Theorem~\ref{thm:HH3.2},
\begin{equation*}
    S'_{n,k_{n}}/\left(\sum_{i=1}^{k_{n}}Y_{n,i}^{2}\right)^{1/2} \stackrel{d}{\to} N(0,1),
\end{equation*}
which, together with \eqref{eq:HH3.6.proof}, implies \eqref{eq:stableCLT_conclusion}.
\end{proof}

\end{appendix}

\bibliographystyle{amsplain}
\bibliography{mybibfile_cleaned} 

@article{christensen2026posterior,
  title={Posterior uncertainty for kernel density estimates},
  author={Christensen, D. and Svardal, T. and R{\o}nneberg, L. and Stoltenberg, E. A.},
  journal={arXiv preprint arXiv:2607.03927},
  year={2026}
}

@article{hilbert2026predictive,
  title={Predictive Inference via Kernel Density Estimates},
  author={Hilbert, T.},
  journal={arXiv preprint arXiv:2605.14008},
  year={2026}
}

@article{fong2026asymptotics,
  title={Asymptotics for a class of parametric martingale posteriors},
  author={Fong, E. and Yiu, A.},
  journal={Biometrika},
  volume={113},
  number={2},
  year={2026},
  publisher={Oxford University Press}
}

@article{mailler2017measure,
  title={Measure-valued {P{\'o}lya} urn processes},
  author={Mailler, C. and Marckert, J.},
  journal={Electronic Journal of Probability},
  volume={22},
  number={26},
  pages={1--33},
  year={2017}
}

@article{sariev2024characterization,
  title={Characterization of exchangeable measure-valued {P{\'o}lya} urn sequences},
  author={Sariev, H. and Savov, M.},
  journal={Electronic Journal of Probability},
  volume={29},
  pages={1--23},
  year={2024},
  publisher={The Institute of Mathematical Statistics and the Bernoulli Society}
}

@article{sariev2025sufficientness,
  title={Sufficientness postulates for measure-valued {P{\'o}lya} urn sequences},
  author={Sariev, H. and Savov, M.},
  journal={Journal of Applied Probability},
  volume={62},
  number={3},
  pages={925--949},
  year={2025},
  publisher={Cambridge University Press}
}

@article{chorbadzhiyska2025some,
  title={Some developments of exchangeable measure-valued {P{\'o}lya} sequences},
  author={Chorbadzhiyska, Y. R. and Sariev, H. and Savov, M.},
  journal={Bayesian Analysis},
  volume={1},
  number={1},
  pages={1--28},
  year={2025},
  publisher={International Society for Bayesian Analysis}
}

@article{bandyopadhyay2022new,
  title={A new approach to {P{\'o}lya} urn schemes and its infinite color generalization},
  author={Bandyopadhyay, A. and Thacker, D.},
  journal={Annals of Applied Probability},
  volume={32},
  number={1},
  pages={46--79},
  year={2022},
  publisher={Institute of Mathematical Statistics}
}

@article{crimaldi2016fluctuation,
  author    = {Crimaldi, I. and Dai Pra, P. and Minelli, I. G.},
  title     = {Fluctuation theorems for synchronization of interacting {P{\'o}lya}’s urns},
  journal   = {Stochastic Processes and their Applications},
  volume    = {126},
  number    = {3},
  pages     = {930--947},
  year      = {2016},
  publisher = {Elsevier}
}

@article{blackwell1973ferguson,
  author    = {Blackwell, D. and MacQueen, J. B.},
  title     = {Ferguson distributions via {P{\'o}lya} urn schemes},
  journal   = {Annals of Statistics},
  volume    = {1},
  number    = {2},
  pages     = {353--355},
  year      = {1973},
  publisher = {Institute of Mathematical Statistics}
}

@article{west1991kernel,
  author    = {West, M.},
  title     = {Kernel density estimation and marginalization consistency},
  journal   = {Biometrika},
  volume    = {78},
  number    = {2},
  pages     = {421--425},
  year      = {1991},
  publisher = {Oxford University Press}
}

@article{fortini2023prediction,
  author    = {Fortini, S. and Petrone, S.},
  title     = {Prediction-based uncertainty quantification for exchangeable sequences},
  journal   = {Philosophical Transactions of the Royal Society A},
  volume    = {381},
  number    = {2247},
  pages     = {20220142},
  year      = {2023},
  publisher = {The Royal Society}
}

@book{rudin1964principles,
  author    = {Rudin, W.},
  title     = {Principles of mathematical analysis},
  year      = {1964},
  publisher = {McGraw-Hill}
}

@article{devroye2018total,
  author  = {Devroye, L. and Mehrabian, A. and Reddad, T.},
  title   = {The total variation distance between high-dimensional {Gaussians} with the same mean},
  journal = {arXiv preprint arXiv:1810.08693},
  year    = {2018}
}

@article{berti2023kernel,
  author    = {Berti, P. and Dreassi, E. and Leisen, F. and Pratelli, L. and Rigo, P.},
  title     = {Kernel based {Dirichlet} sequences},
  journal   = {Bernoulli},
  volume    = {29},
  number    = {2},
  pages     = {1321--1342},
  year      = {2023},
  publisher = {Bernoulli Society for Mathematical Statistics and Probability}
}

@book{silverman2018density,
  author    = {Silverman, B. W.},
  title     = {Density estimation for statistics and data analysis},
  year      = {2018},
  publisher = {Routledge}
}

@book{van1996weak,
  author    = {van der Vaart, A. W. and Wellner, J. A.},
  title     = {Weak convergence and Empirical Processes},
  year      = {1996},
  publisher = {Springer}
}

@article{garelli2024asymptotics,
  author  = {Garelli, S. and Leisen, F. and Pratelli, L. and Rigo, P.},
  title   = {Asymptotics of predictive distributions driven by sample means and variances},
  journal = {Journal of Applied Probability},
  pages   = {1--27},
  year    = {2026}
}

@article{holmes2023statistical,
  author    = {Holmes, C. C. and Walker, S. G.},
  title     = {Statistical inference with exchangeability and martingales},
  journal   = {Philosophical Transactions of the Royal Society A},
  volume    = {381},
  number    = {2247},
  pages     = {20220143},
  year      = {2023},
  publisher = {The Royal Society}
}

@article{berti2004limit,
  author    = {Berti, P. and Pratelli, L. and Rigo, P.},
  title     = {Limit theorems for a class of identically distributed random variables},
  journal   = {Annals of Probability},
  volume    = {32},
  number    = {3},
  pages     = {2029--2052},
  year      = {2004},
  publisher = {JSTOR}
}

@incollection{robbins1971convergence,
  author    = {Robbins, H. and Siegmund, D.},
  editor    = {Rustagi, J. S.},
  title     = {A convergence theorem for non negative almost supermartingales and some applications},
  booktitle = {Optimizing Methods in Statistics},
  pages     = {233--257},
  year      = {1971},
  publisher = {Academic Press}
}

@article{fortini2020quasi,
  author    = {Fortini, S. and Petrone, S.},
  title     = {Quasi-{Bayes} properties of a procedure for sequential learning in mixture models},
  journal   = {Journal of the Royal Statistical Society: Series B},
  volume    = {82},
  number    = {4},
  pages     = {1087--1114},
  year      = {2020},
  publisher = {Royal Statistical Society}
}

@article{fong2021martingale,
  author  = {Fong, E. and Holmes, C. C. and Walker, S. G.},
  title   = {Martingale posterior distributions},
  journal = {Journal of the Royal Statistical Society: Series B},
  volume  = {85},
  number  = {5},
  pages   = {1357--1391},
  year    = {2023}
}

@article{kal88,
  author    = {Kallenberg, O.},
  title     = {Spreading and predictable sampling in exchangeable sequences and processes},
  journal   = {Annals of Probability},
  volume    = {16},
  number    = {2},
  pages     = {508--534},
  year      = {1988},
  publisher = {JSTOR}
}

@incollection{ald85,
  author    = {Aldous, D. J.},
  title     = {Exchangeability and related topics},
  booktitle = {{\'E}cole d'{\'E}t{\'e} de Probabilit{\'e}s de Saint-Flour XIII -- 1983},
  series    = {Lecture Notes in Mathematics},
  volume    = {1117},
  pages     = {1--198},
  year      = {1985},
  publisher = {Springer}
}

@book{bil99,
  author    = {Billingsley, P.},
  title     = {Convergence of Probability Measures},
  year      = {1999},
  publisher = {Wiley}
}

@article{for00,
  author  = {Fortini, S. and Ladelli, L. and Regazzini, E.},
  title   = {Exchangeability, predictive distributions and parametric models},
  journal = {Sankhy{\=a}: Series A},
  volume  = {62},
  pages   = {86--109},
  year    = {2000}
}

@article{ber12,
  author  = {Berti, P. and Pratelli, L. and Rigo, P.},
  title   = {Limit theorems for empirical processes based on dependent data},
  journal = {Electronic Journal of Probability},
  volume  = {17},
  number  = {9},
  pages   = {1--18},
  year    = {2012}
}

@article{ML74,
  author  = {McLeish, D. L.},
  title   = {Dependent central limit theorems and invariance principles},
  journal = {Annals of Probability},
  volume  = {2},
  number  = {4},
  pages   = {620--628},
  year    = {1974}
}

@book{HH80,
  author    = {Hall, P. and Heyde, C. C.},
  title     = {Martingale Limit Theory and Its Applications},
  year      = {1980},
  publisher = {Academic Press}
}

@article{Cr09,
  author  = {Crimaldi, I.},
  title   = {An almost sure conditional convergence result and an application to a generalized {P{\'o}lya} Urn},
  journal = {International Mathematical Forum},
  volume  = {4},
  number  = {23},
  pages   = {1139--1156},
  year    = {2009}
}

@article{bc25,
  author  = {Battiston, M. and Cappello, L.},
  title   = {Bayesian predictive inference beyond martingales},
  journal = {Journal of Machine Learning Research},
 volume  = {27},
  number  = {236},
  pages   = {1--42},
  year    = {2026}
}

@article{glpr26b,
  title={Flexible learning
with adaptive predictive distributions},
  author={Garelli, S. and Leisen, F. and Pratelli, L. and Rigo, P.},
  journal={Bayesian Analysis (to appear)},
  year={2026}
}

@article{bpr13,
  author  = {Berti, P. and Pratelli, L. and Rigo, P.},
  title   = {Exchangeable sequences driven by an absolutely continuous random measure},
  journal = {Annals of Probability},
  volume  = {41},
  number  = {3B},
  pages   = {2090--2102},
  year    = {2013}
}

@article{bdlpr25,
  author  = {Berti, P. and Dreassi, E. and Leisen, F. and Pratelli, L. and Rigo, P.},
  title   = {A Probabilistic View on Predictive Constructions for {Bayesian} Learning},
  journal = {Statistical Science},
  volume  = {40},
  number  = {1},
  pages   = {25--39},
  year    = {2025}
}

@book{DM82,
  author    = {Dellacherie, C. and Meyer, P. A.},
  title     = {Probabilities and Potential {B}: Theory of Martingales},
  series    = {North-Holland Mathematics Studies},
  volume    = {72},
  year      = {1982},
  publisher = {North-Holland}
}

@article{bcl10,
  author  = {Bassetti, F. and Crimaldi, I. and Leisen, F.},
  title   = {Conditionally identically distributed species sampling models},
  journal = {Advances in Applied Probability},
  volume  = {42},
  issue  = {2},
  pages   = {433--459},
  year    = {2010}
}

@article{bw25,
  author  = {Bissiri, P. G. and Walker, S. G.},
  title   = {Bayesian analysis with conditionally identically distributed sequences},
  journal = {Electronic Journal of Statistics},
  volume  = {19},
  pages   = {1609--1632},
  year    = {2025}
}

@article{fps18,
  author  = {Fortini, S. and Petrone, S. and Sporysheva, P.},
  title   = {On a notion of partially conditionally identically distributed sequences},
  journal = {Stochastic Processes and their Applications},
  volume  = {128},
  number = {3},
  pages   = {819--846},
  year    = {2018}
}

@article{fp25,
  author  = {Fortini, S. and Petrone, S.},
  title   = {Exchangeability, prediction and predictive modeling in Bayesian statistics},
  journal = {Statistical Science},
  volume  = {40},
  number  = {1},
  pages   = {40--67},
  year    = {2025}
}

@article{bdpr21,
  author  = {Berti, P. and Dreassi, E. and Pratelli, L. and Rigo, P.},
  title   = {A class of models for {Bayesian} predictive inference},
  journal = {Bernoulli},
  volume  = {27},
  pages   = {702--726},
  year    = {2021}
}

@article{lpr26,
    author = {Leisen, F. and Pratelli, L. and Rigo, P.},
    title = {Weak convergence in probability of predictive distributions},
    journal = {Electronic Communications in Probability},
    volume ={31},
    pages = {1--15},
    year = {2026}
}

@book{DM78,
  author    = {Dellacherie, C. and Meyer, P. A.},
  title     = {Probabilities and Potential},
  series    = {North-Holland Mathematics Studies},
  year      = {1978},
  publisher = {Hermann}
}

@book{hunt66,
  author    = {Hunt, G. A.},
  title     = {Martingales et processus de Markov},
  series    = {Monographies de la Soci{\'e}t{\'e} math{\'e}matique de France},
  number    = {1},
  year      = {1966},
  publisher = {Dunod},
  address   = {Paris},
  language  = {french}
}

@article{blackwell1962merging,
  author    = {Blackwell, D. and Dubins, L.},
  title     = {Merging of Opinions with Increasing Information},
  journal   = {Annals of Mathematical Statistics},
  volume    = {33},
  number    = {3},
  pages     = {882--886},
  year      = {1962},
  publisher = {Institute of Mathematical Statistics},
  doi       = {10.1214/aoms/1177704456}
}

@book{Kec95,
  author    = {Kechris, A. S.},
  title     = {Classical Descriptive Set Theory},
  series    = {Graduate Texts in Mathematics},
  year      = {1995},
  publisher = {Springer},
  isbn      = {9780387943749}
}


\noindent {\large\textbf{Acknowledgments}}

L.C. acknowledges the support of  the Spanish Ministry of Economy and Competitiveness grant PID2022-138268NB-I00, financed by MCIN/AEI/0.13039/501100011033, FSE+MTM2015-67304-P, and FEDER, EU;  the grant Ramon y Cajal 2022 RYC2022-038467-I, financed by\\ MCIN/AEI/10.13039/501100011033 and FSE+, and the Severo Ochoa Programme for Centres of Excellence in R\&D (Barcelona School of Economics CEX2019-000915-S), funded by \\ MCIN/AEI/10.13039/50110001103. M.B. gratefully acknowledges support for this work from the London Mathematical Society, Research in Pairs grant [grant number 42418]. 

\end{document}